\documentclass[12pt]{article}
\usepackage{amssymb}
\usepackage{amsmath}
\usepackage{graphicx}
\usepackage{fullpage}
\usepackage{amsthm}
\numberwithin{equation}{section}

\title{Some remarkable new Plethystic Operators\\ in the Theory of Macdonald Polynomials}

\author{F.~Bergeron, A.~Garsia, E.~Leven\footnote{Supported by NSF grant DGE 1144086.}, and G.~Xin\footnote{Supported by NSFC(11171231).}}

\begin{document}

\maketitle
 
\begin{abstract}

In the 90's a collection of Plethystic operators were introduced in \cite{pl2}, \cite{pl7} and \cite{pl6} to  solve some Representation Theoretical problems arising from the Theory of Macdonald polynomials. This collection was   enriched in the research that led to the results which appeared in \cite{pl5}, \cite{pl8} and \cite{pl9}. However since some of the identities resulting from these efforts  were eventually not needed,   this additional work remained unpublished. As a consequence of very recent publications \cite{pl4}, \cite{pl11}, \cite{pl21}, \cite{pl19}, \cite{pl20}, a truly remarkable expansion of this theory has taken place. However most of this work  has appeared in a language  that is virtually inaccessible to  practitioners of Algebraic Combinatorics. Yet, these developments have led to a variety of new conjectures in \cite{pl3} in the Combinatorics and Symmetric function Theory of Macdonald Polynomials. The present work results from an effort to obtain in an elementary and  accessible manner all the background  necessary to construct the symmetric function side of some of these new conjectures. It turns out that the  above mentioned  unpublished results provide precisely  the tools needed to carry out this project to its completion. 

\end{abstract}

\setcounter{section}{-1}
 
\section{Introduction}

Our main actors in this development are  the operators $D_k$ introduced in \cite{pl6}, whose action 
on a symmetric function $F[X]$ is defined by setting 
\begin{equation}
  D_k F[X]
 \,=\, F[X+\textstyle{M\over z}]\sum_{i\ge 0}(-z)^ie_i[X]\Big|_{z^k}
 \hskip .5in(\hbox{with $M=(1-t)(1-q)$}).
\label{eq:I.1}
\end{equation}
These operators generate an algebra $\mathcal{A}$ of symmetric function operators with remarkable properties. To state  them we need some preliminary observations and definitions. Let us  denote by $\Lambda$ the space of symmetric functions in the infinite alphabet  $X=\{x_1,x_2,x_3,\ldots \}$ and by $\Lambda^{=d}$ the subspace of homogeneous symmetric functions of degree $d$. It is easy to see from \ref{eq:I.1} that if $F[X]\in \Lambda^{=d}$ then $D_kF[X]\in \Lambda^{=d+k}$. Thus $\mathcal{A} $ is clearly a graded algebra. What is surprising is that $\mathcal{A}$ is in fact bi-graded by simply assigning the generators $D_k$ bi-degree $(1,k)$.

To make this more precise consider first \hbox{$\mathcal{D}=\{D_0,D_1,D_2,D_3,\ldots \}$} as an infinite alphabet, and denote by $\mathcal{L}[\mathcal{D}]$ the linear span of words in $\mathcal{D}$. Now, given this bi-grading of the letters of $\mathcal{D}$, every element $\Pi\in \mathcal{L}[\mathcal{D}]$ has a natural decomposition
$$
\Pi\,=\, \sum_{(u,v)}\Pi_{u,v}
$$
where $\Pi_{u,v}$ denotes the portion of $\Pi$ which is a linear combination of words in $\mathcal{D}$ of total bi-degree $(u,v)$. To show that $\mathcal{A}$ is bi-graded it is necessary and sufficient to prove that  $\Pi$, as an operator, acts by zero on $\Lambda$ if and only if all the  $\Pi_{u,v}$ act by zero. This is one of the very first things we will prove about $\mathcal{A}$.

The connection of $\mathcal{A}$ to the above mentioned developments is that it gives a concrete \hbox{realization} of a proper subspace of the Elliptic Hall Algebra studied by Schiffmann and Vasserot in \cite{pl19}, \cite{pl20} and \cite{pl21}. In particular it contains a distinguished family of operators  $\{Q_{u,v}\}$ of bi-degree given by their index that play a central role in the above mentioned conjectures. For a co-prime bi-degree their  construction is so simple that  we need only illustrate it in a special case. 

\vskip -.02in
\hfill $\vcenter{\hbox{\includegraphics[width=1 in]{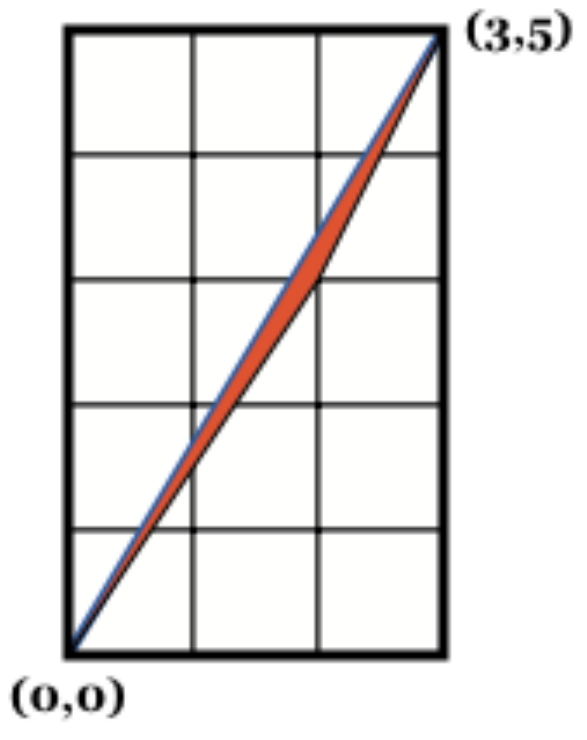}}}$
 
 \vskip -1.23in
\hsize=5.3in
For instance, to obtain $Q_{3,5}$ we start by drawing the $3\times 5$ lattice square with its diagonal (the line  $(0,0)\to (3,5)$, as shown  in the adjacent figure), we then look for the lattice point $   (a,b)$ that is closest to and below the diagonal. In this case $(a,b)=(2,3)$. This yields the decomposition $(3,5)=(2,3)+(1,2)$ and  we set
\vskip -.13in
\begin{equation}
 Q_{3,5}= \textstyle{1\over M} [Q_{1,2}\,,\, Q_{2,3}]
 \,=\, \textstyle{1\over M}\big(
 Q_{1,2} Q_{2,3}-Q_{2,3} Q_{1,2}
 \big).
\label{eq:I.2}
\end{equation}
 
\hsize=6.5in
\hfill 
$\vcenter{\hbox{\includegraphics[width=.75 in]{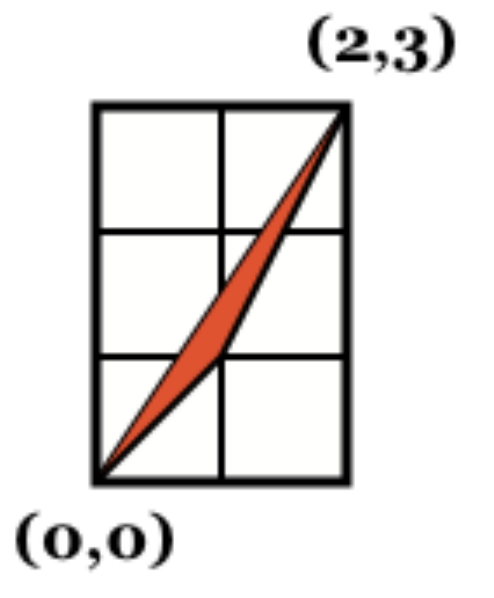}}}$

\vskip -0.9in
\hsize=5.3in
\noindent 
We must next work precisely in the same way with the $2\times 3$ rectangle and, as indicated in the adjacent figure,  obtain the decomposition $(2,3)=(1,1)+(1,2)$
and set 
\begin{equation}
Q_{2,3}\,=\, \textstyle{1\over M} [Q_{1,2}\,,\, Q_{1,1}]
\,=\,
\textstyle{1\over M}\big(
 Q_{1,2} Q_{1,1}-Q_{1,1} Q_{1,2}
 \big).
\label{eq:I.3}
\end{equation}
\hsize=6.5in

\noindent
Now, in this case, we are done, since it turns out that we may set
\begin{equation}
Q_{1,k}\,=\, D_k.
\label{eq:I.4}
\end{equation}
\noindent
In particular by combining \ref{eq:I.2},  \ref{eq:I.3} and \ref{eq:I.4} we obtain 
\begin{align}
Q_{3,5}
&= \textstyle{1\over M^2}
\big(D_2 D_2D_1-2D_2  D_1D_2+D_1D_2 D_2\big).
\label{eq:I.5}
\end{align}

In the general co-prime case $(m,n)$, the precise definition  is based on an elementary number theoretical Lemma  that characterizes the  closest lattice point $(a,b)$ below  the line $(0,0)\to (m,n)$. We then let $(c,d)=(m,n)-(a,b)$ and set
$$
Split(m,n)= (a,b)+(c,d).
$$
This given, we recursively  define
\begin{equation}
Q_{m,n}\,=\, 
\begin{cases}
\textstyle{1\over M}[Q_{c,d},Q_{a,b}] & if $m>1$ and $Split(m,n)=(a,b)+(c,d)$  \cr\cr 
D_n & if $m=1$.
\end{cases}
\label{eq:I.6}
\end{equation}

Our next task is to define the operators $Q_{u,v}$ for any non co-prime pair $(u,v)$. It will be convenient here and after to write such a pair in the form $(km,kn)$ with $(m,n)$ co-prime and $k>1$ the \emph{gcd} of the pair. The problem is that  in this case there are exactly $k$ lattice points, closest to the diagonal of the rectangle $km \times kn$, as we can clearly see in the following display, where we illustrate the case $(m,n)=(3,2)$ and $k=4$. 
$$
\includegraphics[height=.9 in]{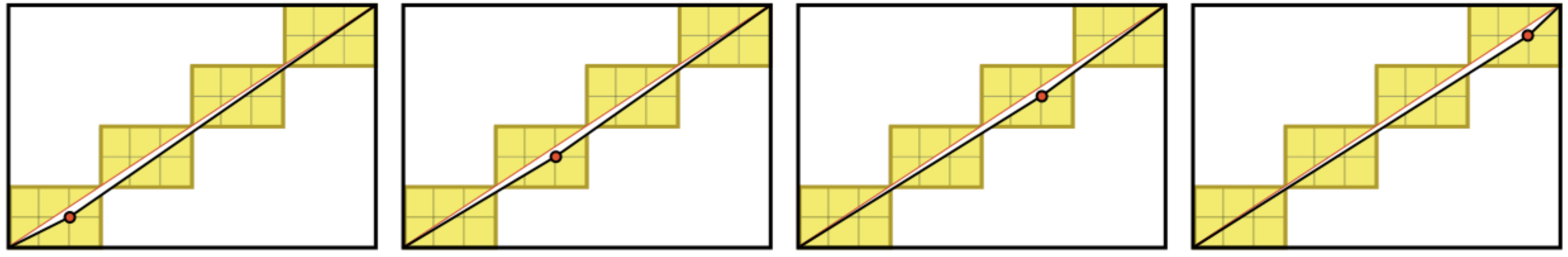}  
$$ 
We see that there are $4$ ways here  to ``\emph{split}'' the  vector $(0,0)\to (4\times 3,4\times 2)$ by choosing a closest lattice point below the diagonal. Namely:
$$
(12,8) = (2,1)+(10,7) =(5,3)+(7,5) = (8,5)+(4,3) = (11,7)+(1,1).
$$

This given, which of the following bracketings should we choose to construct $Q_{4\times 3,4\times 2}$?
$$
[Q_{10,7},Q_{2,1}]\,,\, \hskip 4pt\hskip 4pt
[Q_{7,5},Q_{5,3}]\,,\, \hskip 4pt\hskip 4pt
[Q_{4,3},Q_{8,5}]\,,\, \hskip 4pt\hskip 4pt
[Q_{1,1},Q_{11,7}]. 
$$
The answer is simple: \emph{any one will do},  since all four bracketings give the same operator. This is one of the many identities we need to establish for the operators $Q_{m,n}$. In fact all the pairs $(a,b)$ and $(c,d)$ obtained by splitting
a pair $(km,kn)$, with $(a,b)$ one of the closest lattice points to the segment $(0,0)\to (km,kn)$, are necessarily co-prime.
Our original idea was to prove first the auxiliary identities 
needed to construct the operators $Q_{0,n}$
then obtain all the other needed identities  
as images of the auxiliary identities, under the action of the modular group $G=SL_2[Z]$ on the operators $Q_{m,n}$. In the realization of this plan, the operators $Q_{n,n}$ were more convenient to work with.

More precisely, for a given element 
$\begin{bmatrix} a & c \\  b & d \end{bmatrix}$
we will show that we can set
$$
\begin{bmatrix} a & c \\  b & d \end{bmatrix} Q_{m,n}=Q_{am+cn,bm+dn}.
$$
by proving that two generators of $G$  preserve   all the relations satisfied by the operators $D_k$. 

As we will see,  this is made possible by  means of a very elementary, but surprisingly powerful tool, in Algebraic Combinatorics which has come to be called the \emph{Stanton-Stembridge Symmetrization Trick} (the \emph{SSS Trick} in brief).

By combining the above mentioned auxiliary identities with the action of $G$ we will also be naturally led to the construction of a variety of new additional operators. More precisely there is one operator for each symmetric function $G[X]$, homogeneous of degree $k$ and each co-prime pair $(m,n)$. The resulting operator, which will be  denoted ``${\bf G}_{km,kn},$'' turns out to have  a variety of surprising properties. In fact, computer exploration led to the discovery (in \cite{pl3}) that in many instances the symmetric polynomial ${\bf G}_{km,kn} (-1)^{k(n+1)}$ has a conjectured combinatorial interpretation as an enumerator of certain families of ``rational'' Parking Functions.

One of the most surprising contributions to this branch of Algebraic Combinatorics is a recent deep result \cite{pl18} of Andrei Negut giving a relatively simple but powerful  constant term expression for  the action of the operators $Q_{m,n}$. The reader is referred to the findings concerning the Negut formula that are presented in \cite{pl3} for the reasons we used the word ``powerful'' in this context. Here it has been  one of our priorities to give a straight-forward proof of Negut's formula using only   tools developed in the present treatment of the subject. In our third and final section we present the various results obtained in this  effort. Our main result there is a proof that the  validity  of  the Negut formula  is equivalent to  the statement that a certain quite elementary and completely explicit rational function symmetrizes to zero. It will be seen that this is but  another beautiful  consequence  of the SSS trick. This leads to a computer proof of the  Negut formula in a variety of cases. Moreover, our proof makes it quite clear why and how the so-called ``Shuffle Algebra'' naturally arises in the present context.
In fact, it should be straightforward to extend the machinery used in the proof of the above result to obtain  a proof that, under appropriate definitions, the Shuffle Algebra is isomorphic to the algebra generated by the operators $D_k$. Our presentation terminates with  a proof of the Negut formula  under the specialization at $t=1/q$. \\

\noindent\textbf{ Acknowledgment}

We cannot  overemphasize here the importance of the contribution of Eugene Gorsky and Andrei Negut to the present developments. Without their efforts at translating  their results \cite{pl11}, \cite{pl17}, \cite{pl18}  and results of Schiffmann-Vasserot \cite{pl19}, \cite{pl20}, \cite{pl21} in a language understandable to us, this writing would not have been possible.\\

\section{Notation and Auxiliary identities.}

In dealing with symmetric function identities,  especially those  arising
in the theory of \hbox{Macdonald} Polynomials,  it  is convenient
and often indispensable to use 
plethystic notation. 
This device has a straightforward definition which can be
verbatim implemented in MAPLE  or  MATHEMATICA.
We simply set for
any expression $E=E(t_1,t_2 ,\ldots )$ and any symmetric function $F$ 
\begin{equation}
F[E]= Q_F(p_1,p_2, \ldots )
\Big|_{\textstyle {p_k\to E(\, t_1^k,t_2^k,\ldots )
\rm\hskip 6pt \hskip 4pt for \hskip 4pt all \hskip 4pt k\ge 1}}
\label{eq:1.1}
\end{equation}
where $Q_F$ is the polynomial yielding the expansion of $F$ in 
terms of the power basis. A paradoxical but necessary property of plethystic substitutions is 
 that \ref{eq:1.1} requires  
$ 
p_k[-E]= -p_k[E].
$
This notwithstanding, we will also need to carry out ordinary  changes 
of signs. To distinguish the latter  from the ``\emph{plethystic}''  minus sign,
we will carry out the \hbox{``\emph{ordinary}''} sign change  by
 multiplying  our expressions by  
a new variable ``$\epsilon$'' 
which, outside of the plethystic bracket, is replaced by $-1$. 
Thus we have
$$
p_k[\epsilon E]= \epsilon^k p_k[ E]= (-1)^k p_k[E].
$$
In particular we see that, with this notation, it follows that for any expression $E$ 
and any symmetric function $F$ we may write
\begin{equation}
(\omega F)[E]\,=\, F[-\epsilon E]
\label{eq:1.2}
\end{equation}
where, as customary, ``$\omega$'' denotes the involution that interchanges the 
elementary and homogeneous symmetric function bases.

It will be also good to remind the reader here that many symmetric function identities can be considerably simplified by means of the ``$\Omega$'' notation. For a general expression $E=E(t_1,t_2,\cdots )$
we simply set
$$
\Omega[E]\,=\, exp
\Big(
\sum_{k\ge 1}{p_k[E]\over k}
\Big)
\,=\, 
exp\Big(
\sum_{k\ge 1}{E(t_1^k,t_2^k,\cdots )\over k} 
\Big).
$$
In particular we see that for $X=x_1+x_2+\cdots $
\begin{equation}
\Omega[zX]\,=\, \sum_{m\ge 0}z^m h_m[X]
\label{eq:1.3}
\end{equation}
and  for  $M=(1-t)(1-q)$ we have
\begin{equation}
\Omega[-uM]\,=\, {(1-u)(1-qtu)\over (1-tu)(1-qu)}.
\label{eq:1.4}
\end{equation}

 As in  Macdonald's  \cite{pl16}, for each (french) Ferrers diagram  of a partition $\mu$, and a lattice cell $c\in\mu$ 
 we have  four parameters
$l =l_\mu(c)$, $l'=l'_\mu(c)$, $a =a_\mu(c)$ and  $a'=a'_\mu(c)$ 
called 
\emph{leg, coleg, arm} and \emph{coarm} which  give the 
number of lattice cells of $\mu$  strictly \emph{north, south, east}
and \emph{west}  of $c$.
Denoting by $\mu'$ the conjugate of $\mu$, 
the basic ingredients  we need to keep in mind here are
\begin{align*}
 &
n(\mu)= \sum_{i=1}^{l(\mu)} (i-1)\mu_i\,,\, 
\hskip 4pt\hskip 4pt \hskip 4pt\hskip 4pt
w_\mu(q,t)=\prod_{c\in \mu}(q^{a _\mu(c)} -t^{l _\mu(c)+1})(t^{l _\mu(c)} -q^{a _\mu(c)+1}),
  \cr
 &\hskip 4pt\hskip 4pt\hskip 4pt\hskip 4pt
  T_\mu= t^{n(\mu)}q^{n(\mu')}\,,\,  
\hskip 4pt\hskip 4pt
B_\mu(q,t)=\sum_{c\in \mu}t^{l'_\mu(c)}q^{a'_\mu(c)}
\,,\, \hskip 4pt\hskip 4pt  M=(1-t)(1-q).
\end{align*}
\noindent
Let us recall that the Hall scalar product  is defined by setting
\vskip -.15in
$$
\left \langle p_\lambda\,,\, p_\mu \right \rangle \,=\, 
 z_\mu\hskip 6pt
\chi(\lambda=\mu) \footnote{Here and after we let $\chi({\bf A})=1$ if {\bf A} is true and $\chi({\bf A})=0$ if {\bf A} is false.}
$$
 \noindent
 where $z_\mu$ gives the order of the stabilizer of a permutation with cycle structure $\mu$.

The Macdonald polynomials we work with here are the unique \cite{pl7} symmetric  function basis
$\big\{\widetilde{H}_\mu[X;q,t]\big\}_\mu$ 
 which is  upper triangularly (in dominance order) related
to the modified Schur basis $\{s_\lambda[\textstyle{X\over t-1}]\}_\lambda$ and satisfies the 
orthogonality condition
\begin{equation}
\hskip .5in \left \langle \widetilde{H}_\lambda\,,\, \widetilde{H}_\mu\right \rangle_*\,=\, \chi(\lambda=\mu) w_\mu(q,t),
\label{eq:1.5}
\end{equation}
 \noindent
where $\left \langle\,,\, \right \rangle_*$ denotes a deformation of the Hall scalar product,  defined by setting  
\begin{equation}
\left \langle p_\lambda\,,\, p_\mu \right \rangle_*\,=\, 
(-1)^{|\mu|-l(\mu)}   \prod_i (1-t^{\mu_i})(1-q^{\mu_i})\hskip 6pt z_\mu\hskip 6pt
\chi(\lambda=\mu).
\label{eq:1.6}
\end{equation}
 \noindent 
We will  use here the operator $\nabla$ introduced in \cite{pl1}  by setting 
\begin{equation}
\nabla \widetilde{H}_\mu[X;q,t]\,=\, T_\mu \widetilde{H}_\mu[X;q,t].
\label{eq:1.7}
\end{equation}
We also set for any symmetric function
$F[X]$
\begin{equation}
  D_k^*F[X]
= F[X-\textstyle{\widetilde{M}\over z}]\textstyle  \sum_{i\ge 0}z^ih_i[X]\Big|_{z^k}
\hskip .5in \big(\hbox{with $\widetilde{M}=(1-1/t)(1-1/q)$}\big)
\label{eq:1.8}
\end{equation}
 It will  be convenient to use the symbol  ``$\underline{F}$'' to denote the operator ``\emph{multiplication}'' by a symmetric function $F[X]$. 
 These families of operators were intensively studied in the $90's$ (see  \cite{pl2} and \cite{pl6}) where they gave rise to a variety of conjectures, some of which are  still open. 
 
 In particular it is shown in \cite{pl6} that the operators $D_k$, $D_k^*$,  $\nabla$ 
 and the modified Macdonald polynomials $\widetilde{H}_\mu[X;q,t]$ are related by the following identities.\\

\noindent\textbf{Proposition 1.1} 

 \emph{The operators $D_0,D_0^*$ and $\nabla$  are all self-adjoint
with respect to the $*$-scalar product. Moreover for $k\geq 1$, the operators $\underline{p}_k $, $D_k$ and $D_k^*$ are $*$-scalar product adjoints to $M(-1)^{k-1}p_k^\perp $,  $(-1)^kD_{-k}$ and $(-qt)^kD_{-k}^*$
respectively. We also have}
\begin{equation}
\begin{matrix}
&(i )& \hskip 4pt\hskip 4pt D_0\, \widetilde{H}_\mu= -D_\mu(q,t) \, \widetilde{H}_\mu
\hskip 4pt\hskip 4pt 
&(i )^*&
D_0^*\, \widetilde{H}_\mu= -D_\mu(1/q,1/t) \, \widetilde{H}_\mu
\cr
&(ii)&\hskip 4pt\hskip 4pt  D_k\, \underline{e} - \underline{e}\, D_k = M\, D_{k+1} 
&(ii )^*&
D_k^*\, \underline{e} - \underline{e}\, D_k^* = \, - \widetilde M\, D_{k+1}^* 
\cr
&(iii)&\hskip 4pt\hskip 4pt\nabla \, \underline{e} \nabla^{-1}= -D_1  
&(iii )^*&
\hskip -.06in\nabla \, D_1^* \nabla^{-1}= \underline{e}
\cr
&(iv)&\hskip 4pt\hskip 4pt\nabla^{-1} \, e_1^\perp \nabla =  {\textstyle {1\over M}}D_{-1}   \
&(iv )^*&
\hskip 4pt\hskip 4pt \nabla^{-1} \, D_{-1}^* \nabla = -{\widetilde M\,} e_1^\perp
\end{matrix}
\label{eq:1.9}
\end{equation}
\emph{with $e_1^\perp$ the Hall  scalar product adjoint of multiplication by $e_1$, $ \widetilde M=(1-1/t)(1-1/q)$
 and 
}
\begin{equation}
D_\mu(q,t)= MB_\mu(q,t)-1.
\label{eq:1.10}
\end{equation}
 \noindent
We should mention that recursive applications of \ref{eq:1.9} $(ii)$ and $(ii)^*$ give 
\begin{equation}
a)\hskip 4pt\hskip 4pt D_{k}=
\textstyle{1\over M^k}\sum_{i=0}^k {k \choose r}(-1)^r  \underline{e}_1^r D_0 \underline{e}_1^{k-r}
\,,\,\hskip 4pt\hskip 4pt\hskip 4pt
b)\hskip 4pt\hskip 4pt D_{k}^*=\textstyle{1\over \widetilde M^k} \sum_{i=0}^k {k \choose r}(-1)^{k-r} \underline{e}_1^r D_0^* \underline{e}_1^{k-r}.
\label{eq:1.11}
\end{equation}
 \noindent
For future use, it will be convenient to set
\begin{equation}
a)\hskip 4pt\hskip 4pt \Phi_k= \nabla D_k \nabla^{-1} 
\hskip 4pt\hskip 4pt\hskip 4pt\hskip 4pt\hbox{and}\hskip 4pt\hskip 4pt\hskip 4pt\hskip 4pt
b)\hskip 4pt\hskip 4pt \Psi_k= -(qt)^{1-k}\nabla D_k^* \nabla^{-1}.
\label{eq:1.12}
\end{equation}
This given we have\\

\noindent\textbf{Theorem 1.1}

\emph{The operators $\Phi_k$ and $\Psi_k$ are uniquely determined by the following recursions
\begin{equation}
a)\hskip 4pt\hskip 4pt\Phi_{k+1}=\textstyle{1\over M}[D_1,\Phi_{k}]
\hskip 4pt\hskip 4pt\hskip 4pt\hskip 4pt\hskip 4pt\hskip 4pt\hbox{and}\hskip 4pt\hskip 4pt\hskip 4pt\hskip 4pt\hskip 4pt\hskip 4pt
b)\hskip 4pt\hskip 4pt \Psi_{k+1}=\textstyle{1\over M}[\Psi_{k},D_1]
\label{eq:1.13}
\end{equation}
and initial conditions}
\begin{equation}
a)\hskip 4pt\hskip 4pt\Phi_1=\textstyle{1\over M} [D_1,D_0]
\hskip 4pt\hskip 4pt\hskip 4pt\hskip 4pt\hskip 4pt\hskip 4pt\hbox{and}\hskip 4pt\hskip 4pt\hskip 4pt\hskip 4pt\hskip 4pt\hskip 4pt
b)\hskip 4pt\hskip 4pt\Psi_1=-\underline{e}_1.
\label{eq:1.14}
\end{equation}

\noindent\textbf{Proof}

Note first that, using  \ref{eq:1.9} $(ii)$ and $(iii)$, the definition in 
 \ref{eq:1.12} a) for $k=1$ gives 
 $$
 \Phi_1=\textstyle{1\over M} \nabla\big (D_0e_1-e_1 D_0\big)\nabla^{-1}
 =\textstyle{1\over M}\big(D_0(-D_1)-(-D_1)D_0\big)
 $$
which is another way of writing \ref{eq:1.14} a). The definition in \ref{eq:1.12} b) and   \ref{eq:1.9} $(iii)^*$ give \ref{eq:1.14} b).

Next, conjugating \ref{eq:1.9}  $(ii)$ by $\nabla$ and  using \ref{eq:1.9} $(iii)$ immediately gives \ref{eq:1.13} a). Finally note that, since $\widetilde M =M/qt$
it follows that \ref{eq:1.9}  $(ii)^*$ may be rewritten as 
$$
(qt)^{1-k}D_k^*e_1-e_1(qt)^{1-k}D_k^*\,=\, -M(qt)^{1-k-1}D_{k+1}^*
$$
and \ref{eq:1.13}  b) then follows by conjugating both sides by $\nabla$ and using \ref{eq:1.9}  $(iii)$.

The next identity plays a crucial role in the present development.\\

\noindent\textbf{Theorem 1.2}

\emph{For  $a,b\in Z $ with $a+b>0$ and any symmetric function $F[X]$, we have }
\begin{equation}
\textstyle{1\over M}(D_aD_b^*-D_b^*D_a)\, F[X]\,=\, 
\textstyle{(qt)^b\over qt-1} 
h_{a+b}\big[X(1/qt-1)\big] F[X]. 
\label{eq:1.15}
\end{equation}

A proof of the general identity that includes \ref{eq:1.15} is given in the Appendix.  As  a corollary we obtain\\

\noindent\textbf{Theorem 1.3}

\emph{The operators $\Phi_k$ and $\Psi_k$ defined in \ref{eq:1.12} satisfy the following identity, when $a,b$ are positive integers with sum equal to $n$:}
\begin{equation}
\textstyle{1\over M}[\Psi_b\,,\, \Phi_a]\,=\, \textstyle{ qt \over  qt-1} \hskip 6pt
\nabla \underline{h}_n\big[X(\textstyle{1 \over qt}-1)\big] \nabla^{-1}
\label{eq:1.16}
\end{equation}

\noindent\textbf{Proof}

The identity in \ref{eq:1.15} essentially says that under the given hypotheses the operator $\textstyle{1\over M}(D_b^*D_a-D_aD_b^*) $ acts as multiplication by the symmetric function
$
\textstyle{(qt)^b\over qt-1} 
h_n\big[X(1/qt-1)\big].
$
Thus with our notational conventions \ref{eq:1.15} may be rewritten as
$$
-\textstyle{(qt)^{1-b}\over M}\big(D_b^*D_a-D_a D_b^*\big) \,=\, 
\textstyle{ qt  \over qt-1} 
\underline{h}_n\big[X(1/qt-1)\big].
$$
Conjugating both sides   by $\nabla$ and using 
\ref{eq:1.12} a) and b)  gives \ref{eq:1.16}

Next it is important to keep in mind the following identity which expresses 
 the  action of a sequence of $D_k$ operators  on a symmetric function $F[X]$.\\

\noindent\textbf{Proposition 1.2}
\begin{equation}
D_{a_m}\cdot \cdot \cdot D_{a_1}F[X]=
 F[X+\textstyle{\sum_{i=1}^m{M\over z_i}}]
 \prod_{i=1}^m \Omega[-z_iX]
 {1\over  \prod_{i=1}^m z_i^{a_i}}\,\,
 \hskip -.1in\prod_{1\le i<j\le m}\hskip -.1in\Omega\big[-M\textstyle{z_i/ z_j}\big]
  \Big|_{ z_1^0z_2^0\cdots z_m^0} 
\label{eq:1.17}
\end{equation}

\noindent\textbf{Proof}

It suffices to see what happens when we use \ref{eq:1.1} a) twice:
\begin{align*}
D_{a_2}D_{a_1}F[X]
&=D_{a_2}F[X+\textstyle{M\over z_1}]\Omega[-z_1X]\Big|_{z_1^{a_1}}
= F[X+\textstyle{M\over z_1}+\textstyle{M\over z_2}]\Omega[-z_1(X+\textstyle{M\over z_2})]\Omega[-z_2X]\Big|_{z_1^{a_1}z_2^{a_2}}
\cr
&= F[X+\textstyle{M\over z_1}+\textstyle{M\over z_2}]\Omega[-z_1X]\Omega[-z_2X]
\Omega[-Mz_1/ z_2]
\Big|_{z_1^{a_1}z_2^{a_2}}
\end{align*}

To give a precise general definition of the $Q$ operators  we need 
 the following elementary number theoretical fact that characterizes 
 the  closest lattice point $(a,b)$ below  the line $(0,0)\to (m,n)$. \\

\noindent\textbf{Proposition 1.3}

\emph{For any  pair of co-prime integers $m,n> 1$ there is a unique pair $a,b$ satisfying the following three conditions 
\begin{equation}
(1)\hskip 4pt 1\le a\le m-1 \,,\,\hskip 4pt\hskip 4pt\hskip 4pt\hskip 4pt
(2)\hskip 4pt 1\le b\le n-1 \,,\,\hskip 4pt\hskip 4pt\hskip 4pt\hskip 4pt
(3)\hskip 4pt mb+1=na  
\label{eq:1.18}
\end{equation}
In particular, setting $(c,d)=(m,n)-(a,b)$ we will write for $m,n>1$
\begin{equation}
Split(m,n)\,=\, (a,b)+(c,d)
\label{eq:1.19}
\end{equation}
and otherwise set  
\begin{equation}
a)\hskip 4pt\hskip 4pt Split(1,n)=(1,n-1)+(0,1)\,,\,
\hskip 4pt\hskip 4pt\hskip 4pt\hskip 4pt
b)\hskip 4pt\hskip 4pt Split(m,1)=(1,0)+(m-1,1).
\label{eq:1.20}
\end{equation}
Moreover it follows from our construction that the pairs $(a,b)$ and $(c,d)$ are also co-prime.} \\
\pagebreak
 
\noindent\textbf{Proof } 

When $m,n>1$ the  lattice point that is closest to and strictly below the     diagonal of  the $m\times n$ lattice rectangle
  must  be the unique element of the set
$$
\big\{ (i,\lfloor i\textstyle{n\over m}\rfloor)\, :\, 1\le i \le m-1 \big\}
$$
that minimizes the difference
\begin{equation}
\epsilon_i\,=\, i\textstyle{n\over m}\, - \, \lfloor i\textstyle{n\over m}\rfloor .
\label{eq:1.21}
\end{equation}
In fact, the co-primality of $ m,n$ assures that all these differences are distinct. So the distance minimizer is clearly unique. Next note that if we set
\begin{equation}
k_i= m\epsilon_i= i\, n-m\lfloor i\textstyle{n\over m}\rfloor
\label{eq:1.22}
\end{equation}
then $k_i$ is an integer in the interval
$$
1\le k_i\le m-1
$$
Since all the $k_i$ must be distinct and there are altogether  $m-1$ of them,
 exactly one of them must be equal to $1$. If $k_a=1$ then the minimizing point is $(a,b)$ with $b=\lfloor a \textstyle{n\over m}\rfloor$, and \ref{eq:1.22} for $i=a$ reduces to
\begin{equation}
1\,=\, a\,n\, - \, m \, b.
\label{eq:1.23}
\end{equation}
This proves (1) and (3) of \ref{eq:1.18} and (2) is then an immediate consequence of (1) and
$b=\lfloor a \textstyle{n\over m}\rfloor$. Finally, \ref{eq:1.20}  is simply due to the fact that
in either of these two cases the closest point can be easily identified. 
The co-primality of $(a,b)$ is immediate, since if $(a,b)=(ka',kb')$ for some $k>1$ then $(a',b')$ would be closer to the diagonal
$(0,0)\to(m,n)$. The co-primality of $(c,d)$ holds for the identical reason.

We are now in a position to give  the  definition of the  operators $Q_{m,n}$ that is more suitable for theoretical purposes.\\

\noindent\textbf{Definition 1.1}

\emph{If 
$
Split(m,n)= (a,b)+(c,d) 
$
we recursively set 
\begin{equation}
Q_{m,n}\,=\,\textstyle{1\over M}[Q_{c,d},Q_{a,b}]
\label{eq:1.24}
\end{equation}
with base cases}
\begin{equation}
a)\hskip 4pt\hskip 4pt Q_{1,0}=D_0\hskip 4pt\hskip 4pt\hskip 4pt
\hbox{and}\hskip 4pt\hskip 4pt\hskip 4pt
b)\hskip 4pt\hskip 4pt Q_{0,1}=-\underline{e}_1.
\label{eq:1.25}
\end{equation}

It is easy to see from Proposition 1.3 that recursive applications of \ref{eq:1.24}
will eventually lead to an expression for $Q_{m,n}$ as a polynomial in the non commutative operators $D_0$ and $\underline{e}_1$.

For computer programming purposes the following alternate recursive 
construction is considerably more efficient since, via \ref{eq:1.17}, it gives  all  these operators a plethystic form.\\
\pagebreak

\noindent\textbf{Theorem 1.4} 

\emph{For any pair of co-prime $m,n$ we have}
\begin{equation}
Q_{m,n}\,=\, 
\begin{cases}
\textstyle{1\over M}[Q_{c,d},Q_{a,b}] & if $m>1$
 and $Split(m,n)=(a,b)+(c,d)$\cr 
D_n & if $m=1$.
\end{cases}
\label{eq:1.26}
\end{equation}

\noindent\textbf{Proof}

Since at each application of the Split operation for $m,n>1$ we have both $a\le m-1$ and $c\le n-1$, we will eventually reach the point in the recursion expressed by \ref{eq:1.24} where $m=1$ or $n=1$.
In the first case, \ref{eq:1.20} a) takes over and the identity  in \ref{eq:1.9} $(ii)$ inductively assures that 
$
Q_{1,n}=D_n
$.
In fact,  in the base case  we have, (by \ref{eq:1.9} $(ii)$ for $k=0$)
\begin{equation}
Q_{1,1}\,=\, \textstyle{1\over M}[Q_{0,1},Q_{1,0}]\,=\,  \textstyle{1\over M}[-e_1,D_0]
\,=\,  \textstyle{1\over M}[D_0, e_1]\,=\, D_1
\label{eq:1.27}
\end{equation}
In case $n=1$ and  $m>1$ then   \ref{eq:1.20} b) takes over, yielding
\begin{equation}
Q_{m,1}\,=\, 
 \textstyle{1\over M}[Q_{m-1,1},Q_{1,0}]
 \,=\, \textstyle{1\over M}[Q_{m-1,1},D_0].
\label{eq:1.28}
\end{equation}
Here the base case is reached when $m=2$ yielding
$$
Q_{2,1}\,=\, \textstyle{1\over M}[Q_{1,1},Q_{1,0}]\,=\,  \textstyle{1\over M}[D_1,D_0].
$$

We terminate this section with the following truly surprising and remarkably basic identity in our development.\\

\noindent\textbf{Proposition 1.6}

\emph{For any co-prime pair $m,n$ we have}
\begin{equation}
Q_{m+n,n}\,=\, \nabla Q_{m,n} \nabla^{-1}.
\label{eq:1.29}
\end{equation}

\noindent\textbf{Proof}

We proceed by induction on the size of $min\{m,n\}$. Suppose first that $m,n\ge 2$ and
\begin{equation}
Split(m,n)\,=\, (a,b)+(c,d).
\label{eq:1.30}
\end{equation}
Suppose inductively that we have
\begin{equation}
Q_{a+b,b}=\nabla Q_{a,b}\nabla^{-1}
\hskip 4pt\hskip 4pt\hskip 4pt\hskip 4pt
\hbox{and}
\hskip 4pt\hskip 4pt\hskip 4pt\hskip 4pt
Q_{c+d,d}=\nabla Q_{c,d}\nabla^{-1}.
\label{eq:1.31}
\end{equation}
From \ref{eq:1.30} and Proposition 1.3 it follows that
$$
(1)\hskip 4pt 1\le a\le m-1 \,,\,\hskip 4pt\hskip 4pt\hskip 4pt\hskip 4pt
(2)\hskip 4pt 1\le b\le n-1 \,,\,\hskip 4pt\hskip 4pt\hskip 4pt\hskip 4pt
(3)\hskip 4pt bm+1=na  
$$
Adding $nb$ to both sides of (3) gives
$$
(3')\hskip 4pt\hskip 4pt\hskip 4pt b(m+n)+1=n(a+b),
$$
while  from (1) and (2) it follows that
$$
(1')\hskip 4pt 1\le a+b\le m+n-1 \,,\,\hskip 4pt\hskip 4pt\hskip 4pt\hskip 4pt
(2')\hskip 4pt 1\le b\le n-1 \,,\,\hskip 4pt\hskip 4pt\hskip 4pt\hskip 4pt
$$
But (1'),(2'),(3'), by Proposition 1.3, imply that
$$
Split(m+n,n)\,=\, (a+b,b)+(c+d,d).
$$
This gives
$$
Q_{m+n,n}\,=\, \textstyle{1\over M}[Q_{c+d,d},Q_{a+b,b}]
$$
and from \ref{eq:1.31} we derive that
$$
Q_{m+n,n}\,=\, \nabla Q_{m,n} \nabla^{-1}
$$
completing the induction.

We are left with checking the equality in the cases where 
$m\le 1$ or $n\le 1$. This brings us to the two  identities
\begin{equation}
a)\hskip 4pt\hskip 4pt Split(1,n)=(1,n-1)+(0,1)\,,\,
\hskip 4pt\hskip 4pt\hskip 4pt\hskip 4pt
b)\hskip 4pt\hskip 4pt Split(m,1)=(1,0)+(m-1,1).
\label{eq:1.32}
\end{equation}
The common base case is $(1,1)$. There we must show that
$$
Q_{1,1}\,=\, \nabla Q_{0,1} \nabla ^{-1}.
$$
But by \ref{eq:1.25} b)  and \ref{eq:1.26} this is
$$
D_1\,=\,  -\nabla e_1 \nabla ^{-1}
$$  
which is \ref{eq:1.9}  $(iii)$. We can thus proceed by induction in each case. Now for  case a)  we have
$$
Q_{1,n}\,=\, \textstyle{1\over M}[Q_{0,1},Q_{1,n-1}].
$$
Assuming that the result is true for $n-1$ gives 
$$
\nabla Q_{1,n}\nabla^{-1} \,=\, \textstyle{1\over M}
[\nabla  Q_{0,1}\nabla^{-1},\nabla Q_{1,n-1}\nabla^{-1}]
\,=\, 
 \textstyle{1\over M}
[Q_{1,1},Q_{n,n-1}].
$$
Since 
 $
Split(n+1,n)=(n,n-1)+(1,1)
 $
we  see that
$
 \textstyle{1\over M}
[Q_{1,1},Q_{n,n-1}]
= Q_{n+1,n}.
$
This completes the induction in case a).
Proceeding again by induction in case b) we get
$$
\nabla Q_{m,1}\nabla^{-1} \,=\, \textstyle{1\over M}
[\nabla  Q_{m-1,1}\nabla^{-1},\nabla Q_{1,0}\nabla^{-1}]
\,=\, 
 \textstyle{1\over M}
[Q_{m,1},Q_{1,0}]
$$
and again \ref{eq:1.32} b) gives
 $
\nabla Q_{m,1}\nabla^{-1}=
Q_{m+1,1}
$
completing  the induction and our proof.

\section{The action of $\bf SL_2[Z]$ on the algebra generated by the operators $\bf D_k$.}

To extend the definition of the $Q$ operators 
to non-coprime  pairs of indices we  need to make use 
of the action of $SL_2[Z]$   on the operators $Q_{m,n}$.
More precisely,  
 for
 $\begin{bmatrix} a & c \\ b & d \end{bmatrix} \in SL_2[Z] $
and  any co-prime pair $(m,n)$  we want
\begin{equation}
\begin{bmatrix} a & c \\ b & d \end{bmatrix} Q_{m,n}= Q_{am+cn,bn+dn}.
\label{eq:2.1}
\end{equation}
For this it is sufficient to justify setting
\begin{equation}
NQ_{m,n}=Q_{m+n,n}
\hskip 4pt\hskip 4pt\hskip 4pt\hskip 4pt\hbox{and}\hskip 4pt\hskip 4pt\hskip 4pt\hskip 4pt
SQ_{m,n}= Q_{m,n+m}
\label{eq:2.2}
\end{equation}
for the generators
$$
N=\begin{bmatrix} 1 & 1 \\ 0 & 1 \end{bmatrix}
\hskip 4pt\hskip 4pt\hskip 4pt\hskip 4pt\hbox{and}\hskip 4pt\hskip 4pt\hskip 4pt\hskip 4pt
S=\begin{bmatrix} 1 & 0 \\ 1 & 1 \end{bmatrix}.
$$ 
Since every operator $Q_{m,n}$ is a polynomial in the operators
$D_k=Q_{1,k}$ we will define this action 
on the algebra  generated by the $D_k$ by setting
\begin{equation}
N(D_{k_1}D_{k_2}\cdots D_{k_r})
= Q_{1+k_1,k_1}Q_{1+k_2,k_2}\cdots Q_{1+k_r,k_r}
\label{eq:2.3}
\end{equation}
and
\begin{equation}
S(D_{k_1}D_{k_2}\cdots D_{k_r})
\,=\, Q_{ 1,k_1+1}Q_{ 1,k_2+1}\cdots Q_{ 1,k_r+1}
\label{eq:2.4}
\end{equation}

\noindent
For this action to be well defined it is necessary and sufficient that if any polynomial in the $D_k$ that acts by zero on symmetric functions, then it has an image under $N$ and $S$ which also acts by zero. Now it happens that this fact can be proved by elementary means.

To begin, notice that the identity \ref{eq:1.29}  allows us
to rewrite \ref{eq:2.3} as
\begin{equation}
N(D_{k_1}D_{k_2}\cdots D_{k_r})
=\nabla(  D_{k_1}D_{k_2}\cdots D_{k_r})\nabla^{-1}
\label{eq:2.5}
\end{equation}
which immediately implies the desired property for the action of $N$, since any symmetric function operator that acts by zero has an image under conjugation by $\nabla$ which also acts by zero. 

To prove that $S$ has the desired property we will make use 
of a simple observation which has come to be referred 
to as the \emph{Stanton-Stembridge Symmetrization Trick}.

It may be stated as follows\\

\noindent
\textbf{SSS Trick}

\emph{For  a Laurent polynomial  $F(z_1,z_2,\dots,z_m)$ we have
\begin{equation}
F(z_1,z_2,\dots,z_m)\,  \Big|_{ z_1^0z_2^0\cdots z_n^0} \,=\, 0 
\label{eq:2.6}
\end{equation}
if and only if 
\begin{equation}
Sym_m
F(z_1,z_2,\dots,z_m)\,  \Big|_{ z_1^0z_2^0\cdots z_m^0}\,=\,
0 
\label{eq:2.7}
\end{equation}
where ``$Sym_m$'' is the idempotent that symmetrizes  with 
respect to the variables $z_1,z_2,\dots,z_m$.  
}

It is important to notice  that consequently this is also valid when
 $F$ is a formal power series in other variables with coefficients Laurent polynomials in $z_1,z_2,\dots,z_m$. The surprising circumstance is that quite often the identity in \ref{eq:2.7} turns out to be a consequence of the more encompassing identity
 $$
 Sym_m F(z_1,z_2,\dots,z_m) \,=\, 0
 $$
which is sometimes easier to prove than \ref{eq:2.7}. 
A beautiful example of this type of circumstance is given by the following result which is crucial in our development.

Here and after it will be convenient to use $Z_k$ as an abbreviation for the alphabet $z_1,z_2,\ldots z_k$.\\

\noindent\textbf{Theorem 2.1} 

\emph{Suppose that    
$ 
FP_k(Z_k)  
$  (for each $1\le k\le m$) is  formal power series in other variables with coefficients Laurent polynomials in $z_1,z_2,\ldots ,z_k$. Then, for all symmetric functions $F[X]$ we have
\begin{equation}
\sum_{k=1}^m\Big(
F[X+\textstyle{\sum_{i=1}^k{M\over z_i}}]
\Omega[-Z_kX]
 \,FP_k(Z_k) 
 \hskip -.1in\prod_{1\le i<j\le k}\hskip -.1in\Omega\big[-M\textstyle{z_i/ z_j}\big]\Big)
  \Big|_{ z_1^0z_2^0\cdots z_m^0} \,=\, 0
\label{eq:2.8}
\end{equation}
if and only if
\begin{equation}
Sym_{k}
\Big(
FP(Z_k)\hskip -.1in  \prod_{1\le i<j\le k}
\hskip -.1in \Omega\big[-Mz_i/z_j\big]
\Big)\,=\, 0
\hskip .5in(\hbox{for $1\le k\le m$}).
\label{eq:2.9}
\end{equation}
In particular it follows that 
the operator
$
\mathbf{V}=\sum_{k=1}^m \mathbf{V}_k
$
with 
$
\mathbf{V}_k= \sum_{a}c^{(k)}_{a_1,a_2,\ldots ,a_k}
D_{a_k}\cdots
 D_{a_2}D_{a_1}
$
acts by zero on symmetric polynomials  if and only if, setting
\begin{equation}
\Pi_{\mathbf{V}_k}(Z_k)\,=\, 
\sum_{a}c_{a_1,a_2,\ldots ,a_k}\hskip 6pt 
 \textstyle {1 \over z_1^{a_1} z_2^{a_2}\cdots  z_k^{a_k}},
\label{eq:2.10}
\end{equation}
we have}
\begin{equation}
Sym_{k}
\Big(
\Pi_{\mathbf{V}_k}(Z_k)\hskip -.1in  \prod_{1\le i<j\le k}
\hskip -.1in \Omega\big[-Mz_i/z_j\big]
\Big)\,=\, 0.
\label{eq:2.11}
\end{equation}

\noindent\textbf{Proof}
 
Notice that since for any variable $u$ we have
 $\,
\Omega[-uM]= 1+ \sum_{r\ge 1}u^r h_r[-M]
\, $
it follows that for any $1\le k\le m$ 
$$
 \hskip -.1in\prod_{1\le i<j\le k}\hskip -.1in\Omega\big[-M\textstyle{z_i/ z_j}\big]
$$
is a formal power series in $q,t$ with coefficients Laurent polynomials in $z_1,z_2,\ldots ,z_k$. Likewise  the    expression
\begin{equation}
 F[X+\textstyle{\sum_{i=1}^k{M\over z_i}}]
  \Omega[-Z_kX]\,=\,  F[X+\textstyle{\sum_{i=1}^k{M\over z_i}}]
 \prod_{i=1}^k \Omega[-z_iX]
\label{eq:2.12}
\end{equation}
 may be viewed as a formal power series in the variables in  $X$ and $q,t$ with coefficients Laurent polynomials in $z_1,z_2,\dots ,z_k$. Thus the SSS Trick applies and we can derive from \ref{eq:2.8} and the $S_k$ symmetry 
of the expression in \ref{eq:2.12} that we will have \ref{eq:2.8} if and only if 
\begin{equation}
\sum_{k=1} ^mF[X+\textstyle{\sum_{i=1}^k{M\over z_i}}]
 \  \Omega[-Z_kX]\hskip 6pt
G_k[Z_k]
  \Big|_{ z_1^0z_2^0\cdots z_m^0} = 0
\label{eq:2.13}
\end{equation}
where for convenience we have set
$$
G_k[Z_k]
= Sym_k
 \Big(
 FP_k(Z_k)\,
 \hskip -.1in\prod_{1\le i<j\le k}\hskip -.1in\Omega\big[-M\textstyle{z_i/ z_j}\big]
 \Big).
 $$

In particular we can immediately see that \ref{eq:2.9} implies \ref{eq:2.8}.
We must next  show that the converse is also true. 
To carry this out it will be convenient to set
$$
 \textstyle{1\over z_1}+\textstyle{1\over z_2}+\cdots +\textstyle{1\over z_k} 
\,=\, Z_k^{(-1)}.
$$
This given, note that  the identity
$$
p_a[MZ_k^{(-1)}]= p_a[X+MZ_k^{(-1)}]-p_a[X]
\hskip 4pt\hskip 4pt\hskip 4pt\hskip 4pt (\hbox{for all $a\ge 1$})
$$
gives for $\lambda=(\lambda_1,\lambda_2,	\ldots ,\lambda_l)$
\begin{align*}
p_\lambda[MZ_k^{(-1)}]
&=\prod_{i=1}^l \Big(p_{\lambda_i}[X+MZ_k^{(-1)}]-p_{\lambda_i}[X]\Big) \cr
&=\sum_{S \subseteq \{1,l\}}
\prod_{i\in S} p_{\lambda_i}[X+MZ_k^{(-1)}]
\prod_{i\in \{1,l\}-S}\big( -p_{\lambda _i}[X]\big).
\end{align*}
Using the fact that $p_\lambda[MZ_k^{(-1)}]=p_\lambda[M]p_\lambda[Z_k^{(-1)}]$, we may write
$$
p_\lambda[Z_k^{(-1)}]\,=\,
 \sum_{\mu\preceq \lambda}p_\mu\big[X+MZ_k^{(-1)}\big]c_{\lambda,\mu}[X] 
$$
where ``$\preceq$'' means that all the parts of $\mu$ are parts of $\lambda$,
and, more importantly,  the coefficients $ c_{\lambda,\mu}[X] $ do  not depend  on $k$. Thus a multiple use of  \ref{eq:2.13} with  $F=p_\mu$ for all $\mu\preceq \lambda$ gives
$$
\sum_{k=1} ^m p_\lambda[Z_k^{(-1)}]
 \  \Omega[-Z_kX]\hskip 6pt
G_k[Z_k]
  \Big|_{ z_1^0z_2^0\cdots z_m^0} = 0.
$$
Since $\{p_\lambda[X]\}_\lambda$  is a symmetric function basis it follows from this
that for all symmetric functions $F[X]$ we must also have
$$
\sum_{k=1} ^m F[Z_k^{(-1)}]
 \  \Omega[-Z_kX]\hskip 6pt
G_k[Z_k]
  \Big|_{ z_1^0z_2^0\cdots z_m^0} = 0.
$$

Now notice that for all $k<m$ we have $e_m[Z_k^{(-1)}]=0$. It follows from this that setting $F[X]=e_m^l[X]$ for any $l>1$ but otherwise arbitrary we must have
$$
 \  \Omega[-Z_mX]\hskip 6pt
G_m[Z_m]
 \Big|_{ z_1^l z_2^l \cdots z_m^l} = 0.
$$
Now the expansion $\Omega[-Z_mX]=\sum_{\lambda}m_\lambda[Z_m]h_\lambda[-X]$ together with
the fact that (for $X$ an infinite alphabet) the collection $\{h_\lambda[-X]\}_\lambda$ is a symmetric function basis (thus independent) allows us to conclude that  
for arbitrary $\lambda$ we must have
\begin{equation}
m_\lambda[Z_m]\hskip 6pt G_m[Z_m]
 \Big|_{ z_1^l z_2^l \cdots z_m^l}  = 0.
\label{eq:2.14}
\end{equation}
It follows from our hypotheses that  $G_m[Z_m]$ is a formal power series in other variables with coefficients  Laurent polynomials in $z_1,z_2,\dots,z_m$. Thus if $P[Z_m]$
is any one of these coefficients, from \ref{eq:2.14} we derive that we must also have
 $$
m_\lambda[Z_m]\hskip 6pt P[Z_m]
 \Big|_{ z_1^l z_2^l \cdots z_m^l}  = 0.
$$

We claim that the arbitrariness  of $\lambda$ and $l$ forces the vanishing of  $P[Z_m]$.
To see this note that we may make the substitution
$m_\lambda[Z_m]= \sum_{\lambda(p)=\lambda}z_1^{p_1}z_2^{p_2}\cdots z_m^{p_m}$
where ``$\lambda(p)=\lambda$'' means that the non-zero parts of the weak composition $p$   rearrange to 
the parts of $\lambda$, and obtain
$$
 \sum_{\lambda(p)=\lambda}P[Z_m]\hskip 6pt 
 \Big|_{{ z_1^{l  -p_1 }z_2^{l  -p_2}\cdots z_m^{l  -p_m}}}
 = 0.
$$ 
Since the symmetry of $G_m[Z_m]$ in $z_1,z_2,\dots,z_m$ implies that also $P[Z_m]$
is symmetric in $z_1,z_2,\dots,z_m$, all the above coefficients must be the same. This implies that 
$$
P[Z_m]\hskip 6pt 
 \Big|_{{ z_1^{l  -p_1 }z_2^{l  -p_2}\cdots z_m^{l  -p_m}}}
 = 0
 \hskip .5in(\hbox{for $\lambda(p)=\lambda$}).
$$
But the arbitrariness of $l$ and $\lambda$ gives that we have 
$$
P[Z_m]\hskip 6pt 
 \Big|_{{ z_1^{ q_1 }z_2^{ q_2}\cdots z_m^{q_m}}}
 = 0
 \hskip .5in(\hbox{for all integral vectors $q=(q_1,q_2,\ldots ,q_n)$})
$$
Thus $P[Z_m]$ must identically vanish as asserted. 
Since this holds true for every coefficient of $G_m[Z_m]$ we are led to the conclusion  that  
$$
G_m[Z_m]\,=\, 0.
$$

This not only proves the special case $k=m$ of \ref{eq:2.9} but sets us up for an induction argument on $m$ with base case $m=1$ which is also a particular subcase of the case we have just dealt with. Our proof is thus complete.

As a corollary we obtain\\

\noindent\textbf{Theorem 2.2}

\emph{An operator $\mathbf{V}=\sum_{k=1}^m \mathbf{V}_k$ with 
$$
\mathbf{V}_k=\sum_{a}c^{(k)}_{a_1,a_2,\ldots ,a_k}
D_{a_k}\cdots
 D_{a_2}D_{a_1}
$$
acts by zero on symmetric polynomials  if and only if each of the operators
 $S\mathbf{V}_k$ acts by zero.}\\

\noindent\textbf{Proof}

By Theorem 1.1 $\mathbf{V}$ acts by zero if and only if for each $1\le k\le m$ we have
\begin{equation}
Sym_{k}
\Big(
\Pi_{\mathbf{V}_k}(z_1,z_2,\ldots ,z_k)\hskip -.1in  \prod_{1\le i<j\le k}
\hskip -.1in \Omega\big[-Mz_i/z_j\big]
\Big)\,=\, 0
\label{eq:2.15}
\end{equation}
and $S \mathbf{V}_k$ acts by zero if and only if
\begin{equation}
Sym_{k}
\Big(
\Pi_{S \mathbf{V}_k}(z_1,z_2,\ldots ,z_k)\hskip -.1in  \prod_{1\le i<j\le k}
\hskip -.1in \Omega\big[-Mz_i/z_j\big]
\Big)\,=\, 0.
\label{eq:2.16}
\end{equation}
But from the definition in  \ref{eq:2.4} it follows that
$$
\Pi_{S \mathbf{V}_k}(z_1,z_2,\ldots ,z_k)\,=\,  {\Pi_{ \mathbf{V}}(z_1,z_2,\ldots ,z_k)\over z_1z_2\cdots z_k},
$$
so we see that \ref{eq:2.15} and \ref{eq:2.16} are equivalent identities.

Combining this with the identity in \ref{eq:2.5} and Theorem 1.4 we can now state \\

\noindent\textbf{Theorem 2.3}

\emph{The identities in \ref{eq:2.3} and \ref{eq:2.4}  define an action of the group $G$ on the algebra generated by the operators $D_k$, with the property that for all $\begin{bmatrix} a & c \\ b & d \end{bmatrix} \in G$ we have
$$
\begin{bmatrix} a & c \\ b & d \end{bmatrix}Q_{m,n}= Q_{am+cn,bn+dn}.
$$
In particular this action preserves all the relations satisfied by the operators $Q_{m,n}$ for $(m,n)$ co-prime.}

We have now all we need to define the operators $Q_{km,kn}$. To begin we have the following basic consequence of Theorem 1.3.\\
\pagebreak

\noindent\textbf{Theorem 2.4}

\emph{For any $k\ge 1$ we have 
$Q_{k+1,k}= \Phi_k$ 
and 
$ Q_{k-1,k}= \Psi_k$. In particular, for all pairs $a,b$ of 
positive integers with sum equal to  n  it follows that}
\begin{equation}
\textstyle{1\over M}[Q_{b+1,b}\,,\, Q_{a-1,a}]\,=\, 
\textstyle{ qt \over qt-1} \hskip 6pt
\nabla \underline{h}_n\big[X(\textstyle{1 \over qt}-1)\big] \nabla^{-1}.
\label{eq:2.17}
\end{equation}

\noindent\textbf{Proof }

In view of \ref{eq:1.12} a), the first equality is a special instance of \ref{eq:1.29}. To prove  the second equality, by Theorem 1.3, we only need  show that the operators  $Q_{k-1,k}$ satisfy the same recursions and  base cases as 
the $\Psi_k$ operators. To begin, note that since
$
Split(k,k+1)=(1,1)+(k-1,k)
$ 
it follows that
$$
Q_{k,k+1}\,=\, \textstyle{1\over M}\big[Q_{k-1,k},Q_{1,1} \big]
\,=\, \textstyle{1\over M}\big[Q_{k-1,k},D_1 \big],
$$
which is \ref{eq:1.13} b) for $Q_{k,k+1}$. However the base case is trivial since by definition  $Q_{0,1}=-\underline{e}_1$. The identity in \ref{eq:2.17} is another way of stating \ref{eq:1.16}. This completes our proof.

This proposition has an avalanche of consequences. In  particular, it plays a crucial role in justifying the  definition of the operators $Q_{km,kn}$. The problem, as we mentioned in the introduction, is that in this case, there are $k$ distinct points  
that are closest to the diagonal $(0,0)\to  (km,kn)$.
inside  the $km\times kn$ lattice rectangle. Correspondingly,  if $Split(m,n)=(a,b)+(c,d)$,
we have the following  $k$ ways to split the vector
$(0,0)\to  (km,kn)$:
$$
\big((u-1)m+a,(u-1)n+b\big)\, + \, \big((k-u)m+c,(k-u)n+d\big)
\hskip 4pt\hskip 4pt\hskip 4pt(\hbox{for $1\le u\le k$}).
$$
Theorems \ref{eq:2.3} and \ref{eq:2.4} allow us to overcome this problem
and at the same time prove an important property of the 
$Q_{km,kn}$ operators.\\

\noindent\textbf{Theorem 2.5}

\emph{If $Split(m,n)=(a,b)+(c,d)$ then we may set for  $k> 1$
and any $1\le u\le k$
\begin{equation}
Q_{km,kn}= \textstyle{1\over M}\big[Q_{ (k-u)m+c,(k-u)n+d
}\,,\,
Q_{(u-1)m+a,(u-1)n+b} 
\big].
\label{eq:2.18}
\end{equation}
Moreover,  letting
$
 \Xi= \begin{bmatrix} a & c \\ b & d \end{bmatrix}
\footnote{\hbox{Notice  $\Xi \in SL_2[Z]$ since   (3) of \ref{eq:1.18}
gives  $ad-bc=1$}}$
we   also have
\begin{equation}
a)\hskip 4pt\hskip 4pt Q_{k,k}= \textstyle{ qt \over qt-1} \hskip 6pt
\nabla \underline{h}_k\big[X(\textstyle{1 \over qt}-1)\big] \nabla^{-1} 
\hskip 4pt\hskip 4pt\hskip 4pt\hskip 4pt
\hbox{and}\hskip 4pt\hskip 4pt\hskip 4pt\hskip 4pt
b)\hskip 4pt\hskip 4pt 
 Q_{km,kn}\,=\, \Xi\, Q_{k,k}
\label{eq:2.19}
\end{equation}
In particular it follows that  for any fixed  $(m,n)$ the operators $\big\{Q_{km,kn}\big \}_{k \geq 1}$
form a commuting family.}\\

\noindent\textbf{Proof}

Note first that for $(m,n)=(1,1)$ we have
$
Split(1,1)=(1,0)+(0,1).
$
Thus the right hand side of  \ref{eq:2.18} becomes for any $1\le u \le k$
\begin{equation}
\textstyle{1\over M}\big[Q_{  k-u   , k-u  +1}
\,,\,
Q_{ u , u-1} \big] \,=\, \textstyle{ qt \over qt-1} \hskip 6pt
\nabla \underline{h}_k\big[X(\textstyle{1 \over qt}-1)\big] \nabla^{-1}
\label{eq:2.20}
\end{equation}
where the last equality is another way of writing \ref{eq:2.17}. We thus immediately see that all these assertions are valid for the co-prime pair $(1,1)$, including \ref{eq:2.19} a). 
To deal with  the case of a general co-prime pair $(m,n)$ we notice that a simple calculation gives
$$
\Xi
\Big[
{k-u\atop k-u+1}
\Big]= \Big[
{m(k-u) +c\atop n(k-u)+d} 
\Big]\,,\, 
\hskip 4pt\hskip 4pt\hskip 4pt\hskip 4pt\hskip 4pt\hskip 4pt\hskip 4pt\hskip 4pt
\Xi
\Big[
{u\atop u-1}
\Big]= \Big[
{m(u-1) +a\atop n(u-1) +b } 
\Big].
$$
Thus $\Xi$ maps the  operators occurring on the left hand side of \ref{eq:2.20} 
onto the operators  occurring on the right hand side of \ref{eq:2.18}. Since all these operators are indexed by co-prime pairs,
all the relations they satisfy are preserved by the action of 
the group $G$. In particular 
   the matrix $\Xi$ will map all the equalities resulting from \ref{eq:2.20} into the desired equalities of the right hand sides of \ref{eq:2.18}.  Thus \ref{eq:2.18} well defines the operator $Q_{km,kn}$ and \ref{eq:2.19} b)
necessarily follows. The asserted commutativity follows
just as well, since 
\ref{eq:2.19} a)  shows that the operators $\big\{Q_{k,k}\big \}_{k\ge1}$
form a commuting family, and the identities expressing these commutativities are preserved by $G$. This completes our proof.

An immediate corollary of Theorem 2.5 is  a  recursive construction of  the action of the operators $Q_{km,kn}$ on a symmetric function $F$.\\

\noindent\textbf{Algorithm}

$\hskip 4pt\hskip 4pt$ Given a pair $(km,kn)$ with $(m,n)$  co-prime and $k\ge 1$:
 
 If $km=1$  then $output=D_{n}F$,
 
 \noindent $\hskip 4pt\hskip 4pt\hskip 4pt $ else
 
\textbf{Step 1:}  Pick the first $1\le a\le m$ such that  $1=na-mb$  where $b=\lceil na/m\rceil-1$.

\textbf{Step 2:} Set $(c,d)=(km,kn)-(a,b)$.

\textbf{Step 3:}  $output=\big(Q_{c,d}Q_{a,b}F\, - \, Q_{a,b}Q_{c,d}F\big)/M$.\\

\noindent
Since all these operators lie in the algebra generated by 
the $D_k$, it follows from \ref{eq:1.17} that their action on a symmetric polynomial  may also be given a completely explicit constant term formula.

More precisely, given any pair $(km,kn)$ we can construct a Laurent polynomial $\Pi_{km,kn}[z_1, \dots ,z_{km}]$ such that for every
symmetric polynomial $F[X]$ we have
\begin{equation}
Q_{km,kn}\, F[X]=
 F \big[X+{\sum_{i=1}^{km} \textstyle{M\over z_i}} \big] 
 \prod_{i=1}^{km} \textstyle{\Omega[-z_iX]
 \Pi_{km,kn}[z_1, ... ,z_{km}]} \hskip -.18in
 \displaystyle\prod_{1\le i<j\le km} \hskip -.18in
 \Omega\big[ \hskip -.05in - \hskip -.05in M\textstyle{z_i \over z_j}\big]
  \Big|_{ z_1^0z_2^0\cdots z_{km}^0} \hskip -.1in.
\label{eq:2.21}
\end{equation}
In fact, the above algorithm naturally leads 
to the following result.\\

\noindent\textbf{Proposition 2.1}

\emph{
A family of   Laurent polynomials $\Pi_{km,kn}$
that may be used in \ref{eq:2.21} can be recursively constructed as follows:}

$\hskip 4pt\hskip 4pt$ Given a pair $(km,kn)$ with $(m,n)$  co-prime and $k\ge 1$:
 
 If $km=1$  then set $\Pi_{1,n}={1\over z_1^n}$,
 
else
 
\textbf{Step 1:}  Pick the first $1\le a\le m$ such that  $1=na-mb$ where  $b=\lceil na/m\rceil-1$.

\textbf{Step 2:} let $(c,d)=(km,kn)-(a,b)$.

\textbf{Step 3:} and set  $\Pi_{km,kn}[Z_{1,km}]=\textstyle{1\over M}
\Big(\Pi_{a,b}[Z_{1,a}]\Pi_{c,d}[Z_{a+1,a+c}]
\, - \, 
\Pi_{c,d}[Z_{1,c}]\Pi_{c,d}[Z_{c+1,c+a}]\Big)$.

\noindent
\emph{where for convenience we have set
$Z_{r,s}=\{z_r,z_{r+1},\ldots, z_s\}$}.\\

\noindent\textbf{Proof}

It suffices to show how two such operators compose after they successively act on a symmetric function.
To this end suppose that 
\begin{align}
&a)\hskip 4pt\hskip 4pt\hskip 4pt
\mathbf{V}_A F[X]\,=\,
F\big[X+{\sum_{i=1}^{a}\textstyle{M\over z_i}}\big]
 \prod_{i=1}^{a} \Omega[-z_iX]
 \Pi_{A}[Z_{1,a}]\,
 \hskip -.1in\prod_{1\le i<j\le a}\hskip -.1in\Omega\big[ \hskip -.05in - \hskip -.05in M\textstyle{z_i \over z_j}\big]
  \Big|_{ z_1^0z_2^0\cdots z_{a}^0}
\cr
&c)\hskip 4pt\hskip 4pt\hskip 4pt
\mathbf{V}_C F[X]\,=\,
F\big[X+{\sum_{i=1}^{c}\textstyle{M\over z_i}}\big]
 \prod_{i=1}^{c} \Omega[-z_iX]
 \Pi_{C}[Z_{1,c}]\,
 \hskip -.1in\prod_{1\le i<j\le c}\hskip -.1in\Omega\big[ \hskip -.05in - \hskip -.05in M\textstyle{z_i \over z_j}\big]
  \Big|_{ z_1^0z_2^0\cdots z_{c}^0}
\label{eq:2.22}
\end{align}
Applying $\mathbf{V}_C$ to both sides of \ref{eq:2.22} a) and using \ref{eq:2.22} c) we may write
\begin{align*}
\mathbf{V}_C \mathbf{V}_A\, F[X]
&\,=\,
 F\big[X+{\sum_{i=1}^{a}\textstyle{M\over z_i}}+{\sum_{i=1}^{c}\textstyle{M\over z_{a+i}}}\big]
\prod_{i=1}^{a} \Omega\big[-z_i(X+{\sum_{i=1}^{c}\textstyle{M\over z_{a+i}}})\big]
 \Pi_{A}[Z_{1,a}]\, \Pi_{C}[Z_{a+1,a+c}]\, \cr
 &\hskip .5in\times
 \hskip -.1in\prod_{1\le i<j\le a}\hskip -.1in\Omega\big[ \hskip -.05in - \hskip -.05in M\textstyle{z_i \over z_j}\big]
   \displaystyle \prod_{i=a+1}^{a+c} \Omega[-z_iX] \hskip -.18in
    \prod_{a+1\le i<j\le a+c}\hskip -.18in\Omega\big[ \hskip -.05in - \hskip -.05in M\textstyle{z_i \over z_j}\big]\,  \Big|_{ z_1^0z_2^0\cdots z_{a+c}^0} \cr
&\,=\,
 F\big[X+{\sum_{i=1}^{a+c}\textstyle{M\over z_i}}\big] 
 \prod_{i=1}^{a+c} \Omega[-z_iX] \Pi_{A}[Z_{1,a}]\, \Pi_{C}[Z_{a+1,a+c}]
\hskip -.18in
\prod_{1\le i<j\le a+c}
\hskip -.18in \Omega\big[ \hskip -.05in - \hskip -.05in M\textstyle{z_i \over z_j}\big]\,  
\Big|_{ z_1^0z_2^0\cdots z_{a+c}^0}
\end{align*}
which shows that the Laurent polynomial for
$\mathbf{V}_C \mathbf{V}_A$ may be taken to be $\Pi_{A}[Z_{1,a}]\, \Pi_{C}[Z_{a+1,a+c}]
$.

It is clear, because of the multiplicity of choices of splitting 
a vector $(0,0)\to (km,kn)$, that the  Laurent polynomial 
needed in \ref{eq:2.21} is not unique. However, this non uniqueness 
goes  deeper than it may be suspected, as the following
identity discovered by  Negut  shows.\\

\noindent\textbf{Theorem 2.6}  \cite{pl18}

\emph{For any co-prime pair $(m,n)$ and symmetric function $F[X]$ we have
\begin{equation}
Q_{m,n}F[X]
= 
F[X+\hskip -.03in {\sum_{i=1}^{m}\textstyle{M\over z_i}}] 
 \prod_{i=1}^{m} \Omega[-z_iX] 
 \prod_{i=1}^m\textstyle{1\over z_i^{e_i(m,n)}}
\displaystyle\prod_{i=1}^{m-1}\textstyle{1\over (1-qtz_i/z_{i+1})}
\hskip -.18in
\displaystyle\prod_{1\le i<j\le m}\hskip -.18in
\Omega[-\textstyle{z_i\over z_j}M]
\Big|_{ z_1^0z_2^0\cdots z_{m}^0}
\label{eq:2.23}
\end{equation}
where for convenience we have set}
\begin{equation}
e_i(m,n)=\textstyle{ \lfloor i{n\over m}\rfloor-\lfloor (i-1){n\over m}}\rfloor  
\label{eq:2.24}
\end{equation}
Later in this writing we will present our progress towards providing an elementary proof of this remarkable identity.
Here it is most appropriate to  present some of the consequences of our experimentation with the right hand side of \ref{eq:2.23}. 

The first surprise is that \ref{eq:2.23} is false if $(m,n)$
is replaced by a non co-prime pair. This given, it is best 
to set for any pair of positive integers $(u,v)$ and
symmetric function $F[X]$
\begin{equation}
{\bf N}_{u,v}F[X]= F[X+\hskip -.03in {\sum_{i=1}^{u}\textstyle{M\over z_i}}] 
 \prod_{i=1}^{u} \Omega[-z_iX] 
 \prod_{i=1}^u\textstyle{1\over z_i^{e_i(u,v)}}
\displaystyle \prod_{i=1}^{u-1}\textstyle{1\over (1-qtz_i/z_{i+1})}
\hskip -.15in
\displaystyle \prod_{1\le i<j\le u}\hskip -.15in
\Omega[-\textstyle{z_i\over z_j}M]
\Big|_{ z_1^0z_2^0\cdots z_{u}^0}
\label{eq:2.25}
\end{equation}
and refer to it as the \emph{Negut operator}.

The next surprise is that computer experimentation
led us to formulate the following \hbox{remarkable}\\

\noindent\textbf{Conjecture 2.1}

\emph{For all $k\ge 1$  and all F[X] we have}
\begin{equation}
{\bf N}_{k,k}F[X]\,=\, \nabla \underline{e}_k\nabla ^{-1} F[X].
\label{eq:2.26}
\end{equation}
 
This given, the relation of Negut's  ${\bf N}_{k,k}$ operator to the 
$Q_{km,kn}$ operators
should be given by the following identity.\\

\noindent \textbf{Theorem 2.6}
\begin{equation}
{\bf N}_{k,k}
= (-1)^k \sum_{\lambda \vdash k} 
m_\lambda \big[\textstyle{qt\over qt-1}\big] \big(\textstyle{1-qt\over qt} \big)^{l(\lambda)}
\prod_{i=1}^{l(\lambda)} Q_{\lambda_i,\lambda_i} 
\label{eq:2.27}
\end{equation}

\noindent\textbf{Proof}

Note first that we may write for any two expressions $A,B$
$$
h_k[AB]\,=\, \sum_{\lambda\vdash k }m_\lambda[B] h_\lambda[A].
$$
Letting $A=X({1\over qt}-1)$ and $B={qt\over qt-1}$ gives
$$ 
(-1)^ke_k[X]\,=\, h_k[-X]\,=\, 
\sum_{\lambda\vdash k} 
m_\lambda\big[\textstyle{qt\over qt-1}\big]
h_\lambda\big[X(\textstyle{1 \over qt}-1)\big].
$$
Thus conjugating both sides by $\nabla$ gives
\begin{align}
\nabla \underline{e}_k\nabla^{-1}
&= (-1)^k \sum_{\lambda\vdash k} 
m_\lambda\big[\textstyle{qt\over qt-1}\big]
\nabla \underline{h}_\lambda\big[X(\textstyle{1 \over qt}-1)\big]\nabla ^{-1}.
\label{eq:2.28}
\end{align}
But using \ref{eq:2.19} a) we  easily derive that
$$
\nabla \underline{h}_\lambda\big[X(\textstyle{1 \over qt}-1)\big]\nabla ^{-1}
\,=\, \big(\textstyle{1-qt\over qt} \big)^{\l(\lambda)}
\prod_{i=1}^{\l(\lambda)}Q_{\lambda_i,\lambda_i} 
$$
and we see that, given \ref{eq:2.26}, the identity in \ref{eq:2.27} is simply another way of writing \ref{eq:2.28}.

Thus it would follow from Conjecture 2.1 that the operators ${\bf N}_{k,k}$  are in the algebra generated by the $D_k$ operators. This fact plus a variety of reasons, including experimental evidence, suggested that for the matrix $\Xi$ of Theorem 2.5
we should have $\Xi \, {\bf N}_{k,k}={\bf N}_{km,kn}$. This given, 
applying $\Xi$  to both sides of \ref{eq:2.27} yields  the following extension of \hbox{Conjecture 2.1}.\\

\noindent\textbf{Conjecture 2.2}

\emph{For all co-prime $(m,n)$ and $k\ge 1$ we have}
\begin{equation}
{\bf N}_{km,kn}= (-1)^k \sum_{\lambda\vdash k} 
m_\lambda\big[\textstyle{qt\over qt-1}\big]\big(\textstyle{1-qt\over qt} \big)^{\l(\lambda)}
\prod_{i=1}^{\l(\lambda)}Q_{\lambda_im,\lambda_in}.
\label{eq:2.29}
\end{equation}

The same sequence of steps carried out in the construction of the operator in the right hand side of \ref{eq:2.29},
 can be used to create an infinite family of operators
in the  algebra generated by the $D_k$ operators.
 In fact we need only replace   $e_k$ by any symmetric function of the same degree in the manipulations
carried out in the proof of Theorem 2.6.

To carry this out it is convenient to set for any partition
$\lambda=(\lambda_1,\lambda_2,\cdots ,\lambda_l)$
$$
h_\lambda[X;q,t]\,=\, (\textstyle{qt\over 1-qt})^l \prod_{i=1}^l h_{\lambda_i}[X(1/qt-1)]
$$
and notice  that the collection $\big\{h_\lambda[X;q,t]\big\}_\lambda$ is a symmetric function basis.

This given, we proceed as follows:\\

\noindent\textbf{Definition 2.1}

\emph{Given  any  symmetric function  $G$ that is homogeneous of degree  k  and any co-prime pair  (m,n)}:

\textbf{Step 1:} Construct the expansion
\begin{equation}
G\,=\, \sum_{\lambda\vdash k}c_\lambda(q,t)h_\lambda[X;q,t].
\label{eq:2.30}
\end{equation}

\textbf{Step 2:}
Set
\begin{equation}
{\bf G}_{km,kn} \,=\, \sum_{\lambda\vdash k}c_\lambda(q,t)\prod_{i=1}^{\l(\lambda)}Q_{m\lambda_i,n\lambda_i}.
\label{eq:2.31}
\end{equation}

\noindent\textbf{Remark 2.1}

It is easily seen that  the operator on  the  right hand side of
\ref{eq:2.29} is simply ${\bf G}_{km,kn}$ for $G=e_k$. This immediately
gives rise to a variety of questions. To begin, are there  ways to modify the definition of the Negut  operator ${\bf N}_{u,v}$ to obtain the action of  ${\bf G}_{km,kn}$ for some other choices of $G$. Secondly, we have ${\bf G}_{k,k} = \nabla \, G \, \nabla^{-1}$ whenever $G$ is of degree $k$. Hence the well-known fact that $\nabla e_k $ is Schur positive combined with the fact that ${\bf N}_{k,k}=\nabla e_k\nabla^{-1}$ makes us wonder what cases of Schur positivity may occur for other choices of $G$. Of course it is experimentally well known that  $\pm \nabla s_\lambda$, with an appropriate choice of the sign, is Schur positive. 
We may then ask what bi-graded $S_n$ modules may have Frobenius characteristics given by the  symmetric polynomials resulting from actions of the operators ${\bf G}_{km,kn}$.

It is also conjectured by Haglund et al  \cite{pl13} that a refinement of the polynomial $\nabla e_k$ may also be obtained as an appropriate enumerator of Parking Functions. Using this conjecture Y. Kim in a recent thesis \cite{pl15} shows that for  an infinite variety of 2-row and 2-column partitions  the polynomial $\pm \nabla s_\lambda$ should also be obtained as an enumerator of Parking Functions.
Can  other choices of $G$ lead to similar findings?
It develops that  these questions have some truly surprising 
answers. The reader is referred to  a forthcoming article \cite{pl3}
where  the  ${\bf G}_{km,kn}$ operators, for a variety of choices of the symmetric function $G$ are shown to be closely connected to the combinatorics of the rational
``Parking Functions'' constructed by Hikita in \cite{pl14}.

\section{The Negut operators and the SSS trick}

The problem we deal with in this section is  best understood 
if we start  with an  example. Suppose we 
want to program on the computer the action of the operator
$Q_{5,3}$. Now using \ref{eq:2.21} for $k=1$ and $(m,n)=(5,3)$ we get
\begin{equation}
Q_{5,3}\, F[X]=
 F\big[X+{\sum_{i=1}^{5}\textstyle{M\over z_i}}\big]
 \prod_{i=1}^{5} \Omega[-z_iX]
 \Pi_{5,3}[z_1,z_2, \ldots ,z_{5}]\,
 \hskip -.1in\prod_{1\le i<j\le 5}\hskip -.1in\Omega\big[ \hskip -.05in - \hskip -.05in M\textstyle{z_i \over z_j}\big]
  \Big|_{ z_1^0z_2^0\cdots z_{5}^0} 
\label{eq:3.1}
\end{equation}

\hsize=4.7in

\noindent
where the Laurent polynomial $ \Pi_{5,3}[z_1,z_2,  \ldots ,z_{5}]$ may be obtained by the recursion in Proposition 2.1.  
In this case it is simpler to construct it directly from the binary tree given in the display on the right. The successive splitting 
depicted by this tree immediately gives  

\hsize=6.5in
\vskip -.75in

\hfill $
\vcenter{\hbox{\includegraphics[width=1.8 in]{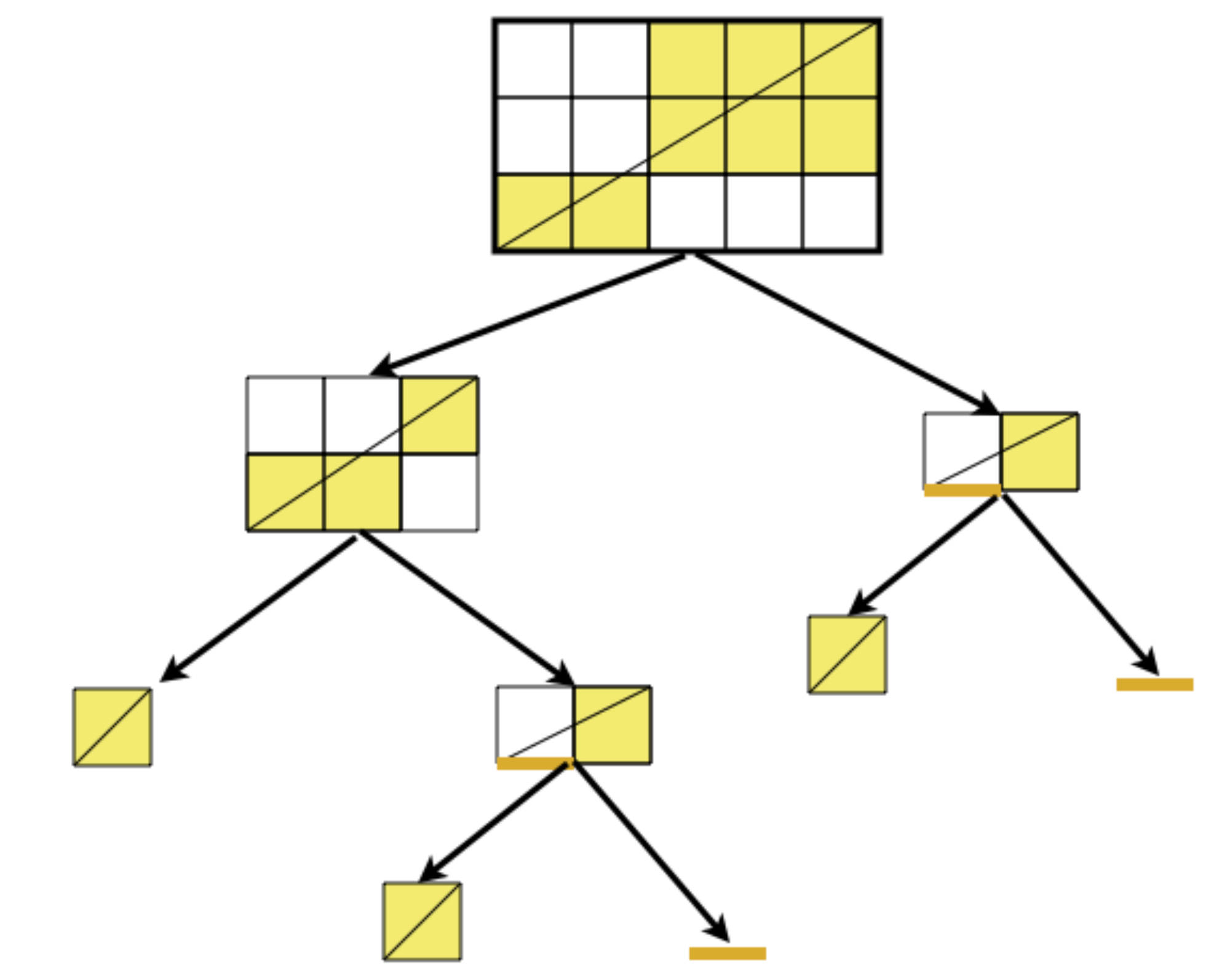}}}
$

\hsize=4.7in
\vskip -.6in
\begin{align*}
Q_{5,3}
&= \textstyle{1\over M} \big[ [D_1,[ D_1,D_0]],[D_1,D_0] \big] \cr 
&= \textstyle{1\over M}
\Big(
\big(
D_1 (D_1D_0-D_0D_1)-(D_1D_0-D_0D_1)D_1  
\big)(D_1D_0-D_0D_1)\, - \,
 \cr
 & \hskip .6in
\, - \, (D_1D_0-D_0D_1)\big(
D_1 (D_1D_0-D_0D_1)-(D_1D_0-D_0D_1)D_1  
\big)
\Big)
\end{align*}

\hsize=6.5in
\noindent
Expanding this out we get
\begin{align*}
Q_{5,3}
=  \textstyle{1\over M^4}&\big(D_1D_1D_0D_1D_0 \, - \,  3 D_1D_0D_1D_1D_0 +2 D_0D_1D_1D_1D_0  \, - \, D_1D_1D_0D_0D_1
\cr
&\hskip 4pt\hskip 4pt\hskip 4pt\hskip 4pt \, + \, 
4 D_1D_0D_1D_0D_1 -3 D_0D_1D_1D_0D_1
\, - \, D_1D_0D_0D_1D_1
\, + \, D_0D_1D_0D_1D_1\big)
\end{align*}
from which we derive that
\begin{align}
\Pi_{5,3}(z_1,z_2,\ldots ,z_5)= 
\textstyle{1\over M^4}
\Big(
{1\over z_2z_4z_5 }
-3{1\over z_2 z_3z_5 }
&+2{1\over  z_2z_3z_4 }
- {1\over  z_1z_4 z_5}
\cr
&
+4{1\over z_1 z_3 z_5}
-3{1\over z_1z_3  z_4}
- {1\over z_1 z_2z_5}
+{1\over  z_1 z_2z_4}
\Big).
\label{eq:3.2}
\end{align}
Now by  \ref{eq:2.25} for $u,v=5,3$, we derive that 
Negut's result is 
\begin{equation}
Q_{5,3} F[X]={\bf N}_{5,3} F[X]
\hskip .5in(\hbox{for all symmetric functions $F[X]$})
\label{eq:3.3}
\end{equation}
with

\begin{equation}
{\bf N}_{5,3}F[X]= F[X+{\sum_{i=1}^{5} \textstyle{M\over z_i}}] 
 \prod_{i=1}^{5} \Omega[-z_iX] 
 \prod_{i=1}^5 \textstyle{1\over z_i^{e_i(m,n)}}
\displaystyle \prod_{i=1}^{4} \textstyle{1\over (1-qtz_i/z_{i+1})}
\hskip -.1in\, \,
\displaystyle \prod_{1\le i<j\le 5}\hskip -.1in\hskip -.1in\,\,
\Omega[-\textstyle{z_i\over z_j}M]
\Big|_{ z_1^0z_2^0\cdots z_{5}^0}
\label{eq:3.4}
\end{equation}

To gauge the simplicity of this formula we need only compute the  monomial $  \prod_{i=1}^5{1\over z_i^{e_i(m,n)}}$.

 \hfill $
\vcenter{\hbox{\includegraphics[width=1.2 in]{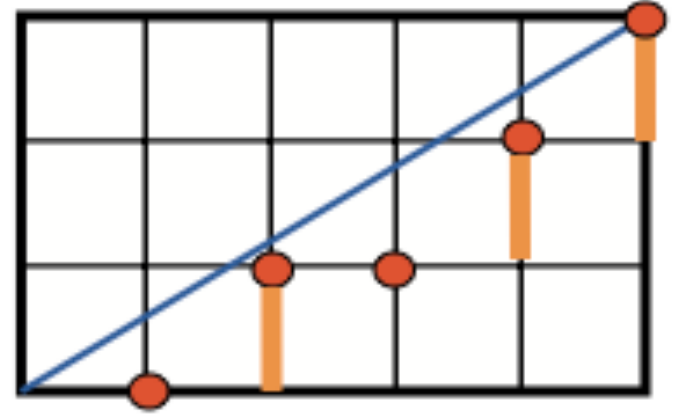}}}
$
 
 \vskip -.85in
 \hsize 5in
 \noindent
 Note that the definition in \ref{eq:2.24}, giving $e_i(m,n)=\textstyle{ \lfloor i{n\over m}\rfloor-\lfloor (i-1){n\over m}}\rfloor$,  geometrically simply means finding, for  each $i$, the highest lattice point $(i,f_i)$ on the line $x=i$ that is below the main diagonal of the lattice rectangle  $m\times n$, then setting $e_i=f_i-f_{i-1}$.
Thus the adjacent display shows that the monomial in \ref{eq:3.3} is simply $z_2z_4z_5$.

 \hsize 6.5in

 Now from Theorem 2.1 we derive that \ref{eq:3.3} can hold true if and only if
 
 \vskip -.2in
\begin{equation}
 Sym_{5}
 \bigg(
 \Big(
 \Pi_{5,3}(z_1,z_2,\ldots ,z_5)-z_2z_4z_5
 \prod_{i=1}^{4}{1\over (1-qtz_i/z_{i+1})}
 \Big)
 \prod_{1\le i<j\le 5}\hskip -.1in\hskip -.1in\,\,
\Omega[-\textstyle{z_i\over z_j}M]
 \bigg)= 0,
\label{eq:3.5}
\end{equation}
 \noindent
an identity that  should be verifiable by computer.

From this example it is easy to deduce the following general result.\\

\noindent\textbf{Theorem 3.1}

\emph{The equality}
\begin{equation}
\hskip 6pt Q_{m,n} F[X]={\bf N}_{m,n} F[X]\hskip 6pt 
\label{eq:3.6}
\end{equation}
\emph{ holds true for all symmetric functions $F[X]$ if and only if}
\begin{equation}
Sym_m
 \bigg( \hskip -.05in \Big(\Pi_{m,n}(z_1,z_2,\dots,z_m)
 \, - \,
\prod_{i=1}^m{1\over z_{ i}^{\lfloor i{n\over m}\rfloor-\lfloor (i-1){n\over m}\rfloor}  }
\prod_{i=1}^{m-1}{1\over (1-qtz_{ i}/z_{ {i+1}})}
\Big) \hskip -.1in
\prod_{1\le i<j\le m} \hskip -.05in
\Omega[-\textstyle{z_{i}\over z_{ j}}M]
\hskip -.02in \bigg)= 0.
\label{eq:3.7}
\end{equation}

 The identity in \ref{eq:3.4} actually was not entirely verifiable on a 
 laptop computer. The problem is not carrying out the symmetrization, but recognizing that the result of symmetrization is actually equal to zero. By setting $z_i=\theta^i$ in  \ref{eq:3.5} then MAPLE 
 is able to recognize that the resulting expression simplifies to zero.
 On the other hand, for the examples in which $m\leq 4$ such as those depicted below
$$
\vcenter{\hbox{\includegraphics[width=3 in]{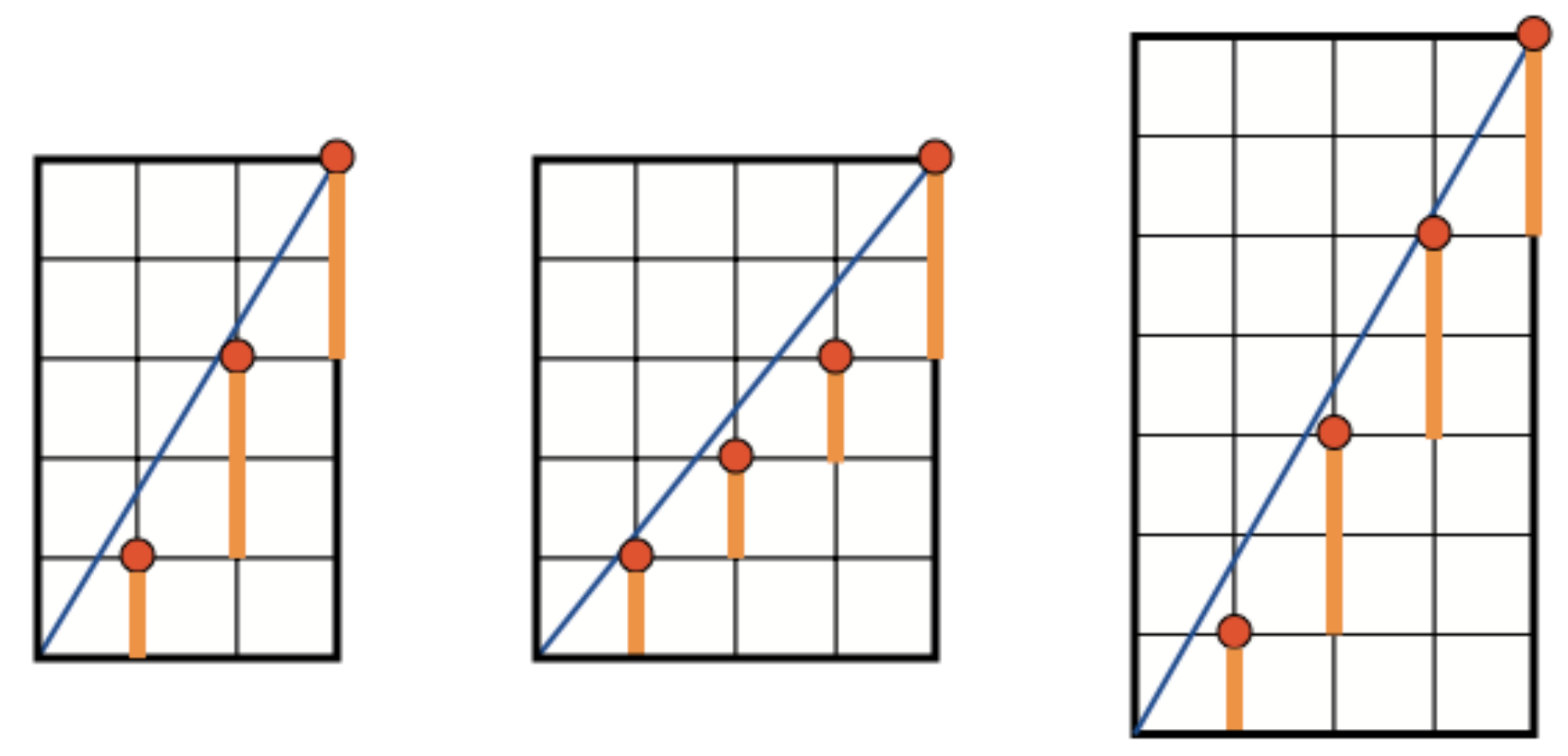}}}
$$ 
the Negut equality $Q_{m,n}={\bf N}_{m,n}$ can be verified even on a laptop in a few seconds. Moreover, as long as $m\le 4$, we  can easily obtain a computer proof of \ref{eq:3.7}.
More precisely\\

\noindent\textbf{Theorem 3.2}

 \emph{For all  co-prime pairs $(m,n)$ with $m\le 4$  we have}
\begin{equation}
 Q_{m,n} ={\bf N}_{m,n}.
\label{eq:3.8}
\end{equation}

\noindent\textbf{Proof}

We might suspect that  $Q_{m,n}-{\bf N}_{m,n}=0$ should 
imply that  $Q_{m,n+m}-{\bf N}_{m,n+m}=0$ by the action  of
the matrix 
$
S=\begin{bmatrix}
 1 & 0 \\
1 & 1 \end{bmatrix}
$.
Of course  we know that $SQ_{m,n}=Q_{m,n+m}$, however at  this moment,  we have no  way to justify acting by $S$
on ${\bf N}_{m,n}$. Nevertheless, the idea can be salvaged,  for  more elementary reason. Simply observe that it follows from \ref{eq:2.4} that
\begin{equation}
\Pi_{m,n+m}(z_1,z_2,\dots,z_m)\,=\, { \Pi_{m,n}(z_1,z_2,\dots,z_m)\over z_1z_2\cdots z_m}.
\label{eq:3.9}
\end{equation}
At the same we can also see that
$$
\prod_{i=1}^m{1\over z_{ i}^{\lfloor i{n+m\over m}\rfloor-\lfloor (i-1){n+m\over m}\rfloor}  }\,=\, 
{1\over z_1z_2\cdots z_m}
\prod_{i=1}^m{1\over z_{ i}^{\lfloor i{n \over m}\rfloor-\lfloor (i-1){n\over m}\rfloor}  }
$$
and we can immediately conclude that the validity of \ref{eq:3.7} for a co-prime pair $(m,n)$ forces the validity of \ref{eq:3.7} for $(m,n+km)$ for any $k\ge 1$. Thus to prove \ref{eq:3.8} for all  pairs $(2,1+2k)$,
 $(3,1+3k)$, $(3,2+3k)$, $(4,1+4k)$, $(4,3+4k)$ it is  sufficient to check it by computer for $k=0$. This can be readily obtained
 in MAPLE or  MATHEMATICA. It is conceivable that by clever means we could succeed in pushing the above computer proof to $m=5$,  but beyond that point it is better to proceed by a more powerful   theoretical approach.

To this end, a moment's reflection should make us plainly see how the  ``\emph{Shuffle Algebra}'' arises  within the present context. In fact, suppose we define as the ``product''
of two symmetric functions $F[Z_a]$, $G[Z_b]$ as the symmetric function $(F\otimes G)[Z_{a+b}]$ defined by setting
\begin{equation}
F[Z_a]\otimes G[Z_{b}]\,=\, Sym_{a+b}\Big(
F[Z_a]G[Z_{a+1,a+b}]\Omega\big[-MZ_a Z_{a+1,a+b}^{-1}\big] 
\Big)
\label{eq:3.10}
\end{equation}
where for an alphabet $Z$ the symbol ``$Z^{-1}$ '' denotes the sum of the inverses  of its letters.
Note that \ref{eq:3.10} can also be rewritten as
\begin{equation}
F[Z_a]\otimes G[Z_{b}]\,=\, {a!b!\over (a+b)!} 
\sum_{\substack{A+B=[a+b] \\ |A|=a,|B|=b}}
F[Z_A]G[Z_B]\Omega\big[-MZ_A Z_B^{-1}\big] 
\Big).
\label{eq:3.11}
\end{equation}
This given, let us set
\begin{equation}
{\bf U}_{m,n}[Z_m]\,=\,
Sym_m \Big(\Pi_{m,n}[Z_m]\prod_{1\le i<j\le m}\Omega[-Mz_i/z_j]\Big).
\label{eq:3.12} 
\end{equation}

Then  from Proposition 2.1 and \ref{eq:3.10} it follows that\\

\noindent\textbf{Proposition 3.1}

\emph{For all co-prime pairs $(m,n)$ with $m>n$ we have}
\begin{equation}
{\bf U}_{m,n}[Z_m]\,=\, 
\textstyle{1\over M} 
\Big(
{\bf U}_{c,d}[Z_c]\otimes {\bf U}_{a,b}[Z_a]\, - \, 
{\bf U}_{a,b}[Z_a]\otimes{\bf U}_{c,d}[Z_c]
\Big).
\label{eq:3.13} 
\end{equation}

\noindent\textbf{Proof}

Note first that
\begin{align*}
Sym_{a+c}&
\Big(\Pi_{c,d}[Z_c]\Pi_{a,b}[Z_{c+1,c+a}]\prod_{1\le i<j\le a+b}\Omega[-Mz_i/z_j]\Big)\,=\,
\cr
& \,=\,
{c!a!\over(c+a)!}
\sum_\tau \tau \bigg(
 Sym_c \Big(\Pi_{c,d}[Z_c]\prod_{1\le i< j\le c}\Omega[-Mz_i/z_j]\Big)
\cr
&\hskip .5in
\times Sym_{c+1,c+a} \Big(\Pi_{a,b}[Z_{c+1,c+a}]\hskip -.2in \prod_{c+1\le i< j\le c+a}
\hskip -.2in \Omega[-Mz_i/z_j]\Big)
\hskip -.2in\prod_{\substack{1\le i\le c \\ c+1\le j\le c+a}}\hskip -.2in\Omega[-Mz_i/z_j]
\bigg)
\cr
&
\,=\,
{c!a!\over(c+a)!}
\sum_{\tau}\tau \bigg({\bf U}_{c,d}[Z_c]{\bf U}_{a,b}[Z_{c+1,c+a}]
\hskip -.2in\prod_{\substack{1\le i\le c \\ c+1\le j\le c+a}}\hskip -.2in
\Omega[-Mz_i/z_j]
\bigg)
\cr(\hbox {by \ref{eq:3.11} })
&\,=\, {\bf U}_{c,d}[Z_c]\otimes{\bf U}_{a,b}[Z_a]
\end{align*}

where the sum is over the left coset representatives $\tau$ of the subgroup 
$S_c\times S_{c+1,c+a} \subseteq S_{a+c}$. It should now be quite clear that the second term in  \ref{eq:3.13} can be obtained in an entirely analogous manner. Thus to complete  our argument, we need only to use the recursion   
\begin{equation}
\Pi_{m,n}[Z_m]\,=\, \textstyle{1\over M}(\Pi_{c,d}[Z_c]\Pi_{a,b}[Z_{c+1,c+a}]\, - \, 
\Pi_{a,b}[Z_a]\Pi_{c,d}[Z_{a+1,a+c}]).
\label{eq:3.14}
\end{equation}

This suggests an inductive approach to the proof of the equality 
$Q_{m,n}={\bf N}_{m,n}$. Name by showing that if
$Split(m,n)=(a,b)+(c,d)$ then
\begin{equation}
{\bf N}_{m,n}\,=\, \textstyle{1\over M}[{\bf N}_{c,d},{\bf N}_{a,b}].
\label{eq:3.15}
\end{equation}
For convenience let  us set
\begin{equation}
\mathbf{V}_{m,n}[Z_m]\,=\,
Sym_m \Big(\Xi_{m,n}[Z_m]\prod_{1\le i<j\le m}\Omega[-Mz_i/z_j]\Big)
\label{eq:3.16} 
\end{equation}
with
\begin{equation}
\Xi_{m,n}[Z_m]= \prod_{i=1}^m{1\over  z_i^{\lfloor i\textstyle{n\over m}\rfloor-\lfloor (i-1)\textstyle{n\over m}\rfloor } }\prod_{i=1}^{m-1}{1\over 1-qt z_i/z_{i+1}}  
={1\over z_m^n}\prod_{i=1}^{m-1}
{ 
(z_i/z_{i+1})^{-\lfloor i\textstyle{n\over m}\rfloor}
\over 1-qt z_i/z_{i+1}
}.
\label{eq:3.17}
\end{equation}

Note that  \ref{eq:2.25} for $(u,v)=(m,n)$ may be written as  
\begin{equation}
{\bf N}_{m,n}F[X]= F[X+\textstyle{\sum_{i=1}^{m}{M\over z_i}}] 
 \prod_{i=1}^{m} \Omega[-z_iX] 
\hskip 6pt \Xi_{m,n}[Z_m]\hskip -.1in\, \,
\prod_{1\le i<j\le m}\hskip -.1in\hskip -.1in\,\,
\Omega[-\textstyle{z_i\over z_j}M]
\Big|_{ z_1^0z_2^0\cdots z_{m}^0}
\label{eq:3.18}
\end{equation}
and since  $\Xi_{1,n}(Z_1)={1\over z_1}$ it follows that 
${\bf N}_{1,n}=D_1=Q_{1,n}$ for all $n\ge 1$. Thus a proof of \ref{eq:3.15}
is all that is needed to prove the Negut  equality $Q_{m,n}={\bf N}_{m,n}$.

As further evidence of the isomorphism between the Shuffle Algebra and the algebra generated by the $D_k$ operators, we must point out that from Theorems \ref{eq:2.1} and  \ref{eq:3.18} we may easily derive the following.\\ 

\noindent\textbf{Proposition 3.2}

\emph{For a co-prime pair $(m,n)$ with $Split(m,n)=(a,b)+(c,d)$ we have }
$$
{\bf N}_{m,n}= \textstyle{1\over M}[{\bf N}_{c,d},{\bf N}_{a,b}]
\hskip 4pt\hskip 4pt\hskip 4pt  \Longleftrightarrow \hskip 4pt\hskip 4pt\hskip 4pt 
\mathbf{V}_{m,n}[Z_m]= \textstyle{1\over M}
\Big(
\mathbf{V}_{c,d}[Z_c]\otimes \mathbf{V}_{a,b}[Z_a]\, - \, 
\mathbf{V}_{a,b}[Z_a]\otimes\mathbf{V}_{c,d}[Z_c]
\Big)
$$

The following general result enabled us to obtain a
computer proof of Negut's equality \ref{eq:3.6} 
for all $m\le 7$.\\

\noindent\textbf{Theorem 3.3}

\emph{If  $m>n$ is a co-prime pair with $Split(m,n)=(a,b)+(c,d)$ then the identity 
\begin{equation}
{\bf N}_{m,n}= \textstyle{1\over M}[{\bf N}_{c,d},{\bf N}_{a,b}] 
\label{eq:3.19}
\end{equation}
holds true if and only if 
\begin{equation}
Sym_m
\bigg(
{\Xi_{m,n}[Z_m]
\over z_{a+1}z_c}
\Big(
z_{a+1} z_{c+1}
-tz_{a+1}z_c-qz_{a+1}z_c
\, + \,
qtz_az_c
\Big)
\prod_{1\le i<j\le m}\hskip -.1in\hskip -.1in\,\,\Omega[-\textstyle{z_i\over z_j}M]
\bigg)\,=\, 0.
\label{eq:3.20}
\end{equation}
In particular the validity of \ref{eq:3.20} forces the equality}
$$
{\bf N}_{m,n}\,=\, Q_{m,n}.
$$

\noindent\textbf{Proof}

From \ref{eq:3.18} and Theorem 2.1 we derive that \ref{eq:3.19} holds true if and only if
\begin{equation}
\hskip -.02in Sym_m \hskip -.05in \left( \hskip -.05in
\Big( \hskip -.03in
M \Xi_{m,n}[Z_m] \hskip -.03in - \hskip -.03in \big( \Xi_{a,b}[Z_a]\, \Xi_{c,d}[Z_{a+1,a+c}] \hskip -.02in - \hskip -.02in \Xi_{c,d}[Z_c]\Xi_{a,b}[Z_{c+1,c+a}] \big)
\Big) \hskip -.2in
\prod_{1\le i<j\le m} \hskip -.18in \Omega[-\textstyle{z_i\over z_j}M] \hskip -.03in
\right) \hskip -.05in
= 0. \hskip -.18in
\label{eq:3.21}
\end{equation}
To compute the first factor within $Sym_m$ we need to  consider the two cases $a<c$ and $a>c$ which are schematically depicted in the display below.
$$
\vcenter{\hbox{\includegraphics[height=1in,width=5 in]{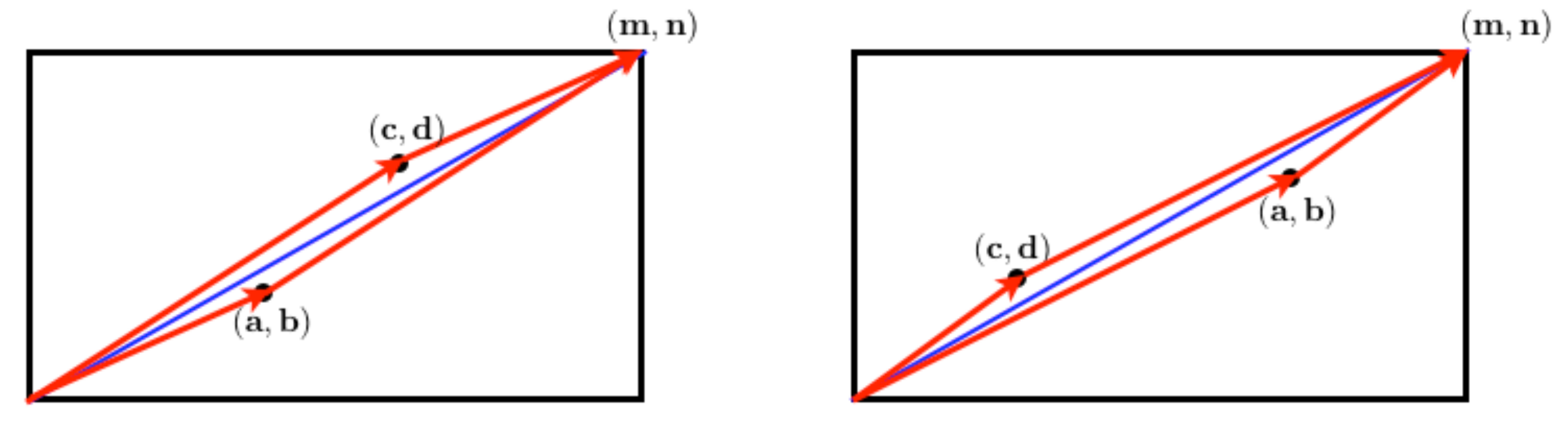}}}
$$
Note first that in each case there are no lattice points within the red parallelogram. Thus the set of highest lattice points $(i,f_i)$ below the  diagonal $(0,0)\to (m,n)$ is the same as the set of highest lattice points $(i,f_i)$ below the  vector sum $(a,b)+(c,d)$, except for $(a,b)$.
This gives
$$
\lfloor i\textstyle{n\over m}\rfloor= \lfloor i\textstyle{b\over a}\rfloor
\hskip 4pt\hskip 4pt  ( \hbox{for $1\le i\le a$})
\hskip 4pt\hskip 4pt\hskip 4pt\hskip 4pt\hskip 4pt\hskip 4pt
\hbox{and}\hskip 4pt\hskip 4pt\hskip 4pt\hskip 4pt\hskip 4pt\hskip 4pt
\lfloor i\textstyle{n\over m}\rfloor= b+\lfloor (i-a)\textstyle{d\over c}\rfloor
\hskip 4pt\hskip 4pt  (\hbox{for $a+1\le i\le a+c=m$}).
$$
Thus \ref{eq:3.17} gives
\begin{align*}
 \Xi_{a,b}[Z_a]\, \Xi_{c,d}[Z_{a+1,a+c}]
 &=
 {1\over z_a^b}\prod_{i=1}^{a-1}
{ 
(z_i/z_{i+1})^{-\lfloor i\textstyle{b\over a}\rfloor}
\over 1-qt z_i/z_{i+1}
}
{1\over z_{a+c}^d}\prod_{i=a+1}^{a+c-1}
{ 
(z_i/z_{i+1})^{-\lfloor( i-a)\textstyle{d\over c}\rfloor}
\over 1-qt z_i/z_{i+1}
}  
\cr
&=
\big({z_a\over z_{a+1}}\big)^b
\big(1-qt{z_a\over z_{a+1}}\big) {1\over z_a^b}
\prod_{i=1}^{a}
{ 
(z_i/z_{i+1})^{-\lfloor i\textstyle{n\over m}\rfloor}
\over 1-qt z_i/z_{i+1}
}
{1\over z_{m}^d}
\prod_{i=a+1}^{m-1}
{ 
(z_i/z_{i+1})^{b-\lfloor i\textstyle{n\over m}\rfloor}
\over 1-qt z_i/z_{i+1}
}  
\cr
&=
{1\over z_{a+1}^b}
\big(1-qt{z_a\over z_{a+1}}\big) {z_m^n\over z_m^d}
\hskip 6pt\Xi_{m,n}[Z_m]
(z_{a+1}/z_m)^b
\,=\, 
\big(1-qt{z_a\over z_{a+1}}\big)
\, \Xi_{m,n}[Z_m].
\end{align*}

For the same reason, in each case 
 the set of highest lattice points $(i,f_i)$ below
the diagonal $(0,0)\to (m,n)$ is the same as the set highest lattice points $(i,f_i)$ below the reversed  vector sum $(c,d)+(a,b)$, except for $(c,d)$.
This gives
$$
\lfloor i\textstyle{n\over m}\rfloor= \lfloor i\textstyle{d\over c}\rfloor -\chi(i=d)
\hskip 4pt  ( \hbox{for $1\le i\le c$})
\hskip 8pt \hbox{and} \hskip 8pt
\lfloor i\textstyle{n\over m}\rfloor= d+\lfloor (i-c)\textstyle{b\over a}\rfloor
\hskip 4pt  (\hbox{for $c+1\le i\le c+a=m$})
$$
thus  again from \ref{eq:3.17} we get
 \begin{align*}
 \Xi_{c,d}[Z_c]\, \Xi_{a,b}[Z_{c+1,c+a}]
 &=
 {1\over z_c^d}\prod_{i=1}^{c-1}
{ 
(z_i/z_{i+1})^{-\lfloor i\textstyle{d\over c}\rfloor}
\over 1-qt z_i/z_{i+1}
}
{1\over z_{c+a}^b}
\prod_{i=c+1}^{ c+a-1}
{ 
(z_i/z_{i+1})^{-\lfloor( i-c)\textstyle{b\over a}\rfloor}
\over 1-qt z_i/z_{i+1}
}  
\cr
&=
\big({z_c\over z_{c+1}}\big)^{d-1}
\big(1-qt{z_c\over z_{c+1}}\big) {1\over z_c^d}
\prod_{i=1}^{c}
{ 
(z_i/z_{i+1})^{-\lfloor i\textstyle{n\over m}\rfloor}
\over 1-qt z_i/z_{i+1}
}
{1\over z_{m}^b}
\prod_{i=c+1}^{m-1}
{ 
(z_i/z_{i+1})^{d-\lfloor i\textstyle{n\over m}\rfloor}
\over 1-qt z_i/z_{i+1}
}  
\cr
&=
{z_c^{-1}\over z_{c+1}^{d-1}}
\big(1-qt{z_c\over z_{c+1}}\big) 
{z_m^n\over z_m^b}
\hskip 6pt \Xi_{m,n}[Z_m]
(z_{c+1}/z_m)^d \cr
&=
z_{c+1}z_c^{-1} \big(1-qt{z_c\over z_{c+1}}\big)
\Xi_{m,n}[Z_m].
\end{align*}

Combining these two identities we get   
\begin{align*}
\Big(
M\, \Xi_{m,n}[Z_m]& \, - \, \big( \Xi_{a,b}[Z_a]\, \Xi_{c,d}[Z_{a+1,a+c}]\, - \, \Xi_{c,d}[Z_c]\Xi_{a,b}[Z_{c+1,c+a}]\big) 
\Big)=
\cr
&=
\Xi_{m,n}[Z_m]
\Big(
(1-t)(1-q)
\, - \,
\big(1-qt\textstyle{z_a\over z_{a+1}}\big)
\, + \, z_{c+1} z_c^{-1} \big(1-qt\textstyle{z_c\over z_{c+1}}\big)
\Big)
\cr
&=
{\Xi_{m,n}[Z_m]
\over z_{a+1}z_c}
\Big(
z_{a+1}z_c(1-t-q+qt) -
\big(z_{a+1}z_c-qt{z_az_c}\big)
+ z_{a+1} z_{c+1} \big(1-qt \textstyle{z_c\over z_{c+1}}\big)
\Big)
\cr
&=
{\Xi_{m,n}[Z_m]
\over z_{a+1}z_c}
\Big(
z_{a+1}z_c(-t-q)
\, + \,
qtz_az_c
\, + \, z_{a+1} z_{c+1}
\Big).
\end{align*}

This shows that \ref{eq:3.21} is equivalent to \ref{eq:3.20} and completes our proof.

To this date we have not yet been able to prove \ref{eq:3.20}  in full generality. However, Theorem 3.3 has the following immediate corollary.\\

\noindent\textbf{Theorem 3.4}

\emph{For any co-prime pair $(m,n)$ we have}
\begin{equation}
{\bf N}_{m,n}\Big|_{t=1/q}\,=\, Q_{m,n}\Big|_{t=1/q}.
\label{eq:3.22}
\end{equation}

\noindent\textbf{Proof}

It is sufficient  to verify \ref{eq:3.20} for $t=1/q$. To begin notice that
$$
\Omega[-uM]\,=\, {q(1-u)^2\over (q-u)(1-qu)}.
$$
This gives 
$$
\prod_{1\le i<j\le m}\hskip -.1in
\Omega[-Mz_i/z_j]\Big|_{t=1/q} = q^{m\choose 2}
\prod_{1\le i<j\le m}\hskip -.1in
{(z_i-z_j)^2\over (qz_j-z_i)(z_j-qz_i)}
$$
which is easily seen to be a symmetric rational function.
Thus we only need  to show that
$$
Sym_m
\Big(
{\Xi_{m,n}[Z_m]\Big|_{t=1/q}
 }
\big(
z_{c+1}/z_c
-(q+1/q) 
\, + \,
z_a/z_{a+1}
\big)
\Big)\,=\, 0
$$
or equivalently that
$$
Asym_m
\Big(
{\Xi_{m,n}[Z_m]\Big|_{t=1/q}\Delta[Z_m]
 }
\big(
z_{c+1}/z_c
-(q+1/q) 
\, + \,
z_a/z_{a+1}
\big)
\Big)\,=\, 0
$$
with $\Delta[Z_m]=\prod_{1\le i<j\le m}(z_j-z_i)$ and ``$Asym_m$'' denoting 
$S_m$
antisymmetrization.

Next notice that since we need only prove \ref{eq:3.22} for 
$1\le n\le m-1$ we can assume, here and after,  that $\textstyle{n\over m}<1$  and we may write
\begin{align*}
\Xi_{m,n}[Z_m]\Big|_{t=1/q}
\Delta[Z_m]
&=
 \prod_{i=1}^m{1\over  z_i^{\lfloor i\textstyle{n\over m}\rfloor-\lfloor (i-1)\textstyle{n\over m}\rfloor } }\prod_{i=1}^{m-1}
 {z_{i+1}\over z_{i+1}- z_i}  \Delta[Z_m]
 \cr
 &=
  \prod_{i=2}^m
 {z_i\over  z_i^{\lfloor i\textstyle{n\over m}\rfloor-\lfloor (i-1)\textstyle{n\over m}\rfloor } }
 \prod_{i=1}^{m-2} \prod_{j=i+2}^{m}(z_j-z_i).
\end{align*}
Thus we are reduced to showing that
$$
Asym_m
\bigg(
\Big( {z_{c+1}\over z_c}
-(q+1/q) 
\, + \,
{z_a\over z_{a+1}}
\Big)  \prod_{i=2}^m
 {z_i\over  z_i^{\lfloor i\textstyle{n\over m}\rfloor-\lfloor (i-1)\textstyle{n\over m}\rfloor } }
 \prod_{i=1}^{m-2} \prod_{j=i+2}^{m}(z_j-z_i)
\bigg)\,=\, 0.
$$

 For the pairs $(m,n)$ with $1\le n\le m-2$ we can prove that the expression to be anti-symmetrized is actually a homogeneous polynomial of degree less than $m \choose 2$, which is the minimum possible degree 
for which $Asym_m$ can yield something other than zero.
For the pair $(m,m-1)$ we will need to prove \ref{eq:3.2}3   by a direct brute force argument. 

Notice first  that for any $1\le i\le m$ we have
$$
\lfloor i\textstyle{n\over m}\rfloor-\lfloor (i-1)\textstyle{n\over m}\rfloor \le 1.
$$
In fact setting $\lfloor i\textstyle{n\over m}\rfloor=r$ and  letting
$\epsilon= i\textstyle{n\over m}-r$ we derive
$$
\lfloor (i-1)\textstyle{n\over m}\rfloor 
\,=\, \lfloor r+\epsilon-\textstyle{n\over m}\rfloor\,=\,
\begin{cases}
r+\lfloor\epsilon-\textstyle{n\over m}\rfloor=r & \hbox{if } \epsilon\ge\textstyle{n\over m}
\cr
r-\lceil \textstyle{n\over m}-\epsilon\rceil=r-1 & \hbox{if } \epsilon<\textstyle{n\over m}.
\end{cases}
$$
This proves \ref{eq:3.2}4.  

Let us now suppose  that $n\le m-2$. 

In view of \ref{eq:3.2}4 to show that the expression inside  $Asym_m$ is a polynomial we need only show that
$$
1)\hskip 4pt\hskip 4pt\lfloor c\textstyle{n\over m}\rfloor-\lfloor (c-1)\textstyle{n\over m}\rfloor 
=0
\hskip 4pt\hskip 4pt\hskip 4pt\hskip 4pt\hskip 4pt\hskip 4pt\hbox{and}\hskip 4pt\hskip 4pt\hskip 4pt\hskip 4pt\hskip 4pt\hskip 4pt
2)\hskip 4pt\hskip 4pt\lfloor (a+1)\textstyle{n\over m}\rfloor-\lfloor a\textstyle{n\over m}\rfloor 
=0
$$
This given, notice that since we clearly have
$\sum_{i=2}^m(\lfloor i\textstyle{n\over m}\rfloor-\lfloor (i-1)\textstyle{n\over m}\rfloor )=n$, the degree of the resulting polynomial must be 
$$
m-1-\sum_{i=2}^m\big(\lfloor i\textstyle{n\over m}\rfloor-\lfloor (i-1)\textstyle{n\over m}\rfloor \big)\, + \, {m-1\choose 2}
\,=\, m-1-n+ {m-1\choose 2}<  {m\choose 2}
$$
as desired to show this polynomial to anti-symmetrize to zero.

To prove \ref{eq:3.2}5 we note that the equality in \ref{eq:1.18} (3),
namely 
$$
na=b\, m+1,
$$
together with $(m,n)=(a,b)+(c,d)$, gives that
$$
nc=md-1.
$$
Thus since ${n+1\over m}<1$
$$
\lfloor c\textstyle{n\over m}\rfloor-\lfloor (c-1)\textstyle{n\over m}\rfloor 
\,=\,
\lfloor d-\textstyle{1\over m}\rfloor-\lfloor d-\textstyle{n+1\over m}\rfloor
\,=\, (d-1)-(d-1)=0.
$$
Likewise  \ref{eq:3.2}6 gives
$$
\lfloor (a+1)\textstyle{n\over m}\rfloor-\lfloor a\textstyle{n\over m}\rfloor
\,=\,
 \lfloor b+  \textstyle{n+1\over m}\rfloor
 -\lfloor b+\textstyle{1\over m}\rfloor\,=\,  0
$$
again since ${n+1\over m}<1$, completing our proof of \ref{eq:3.2}5.

Finally suppose that $n=m-1$. Since in this case
$ Split((m,m-1)=(m-1,m-2)+(1,1)$ and
$$
\prod_{i=2}^m
{z_i\over  z_i^{\lfloor i\textstyle{n\over m}\rfloor-\lfloor (i-1)\textstyle{n\over m}\rfloor }}\,=\, 1
$$
then \ref{eq:3.2}3 reduces to
$$
Asym_m
\bigg(
\Big( {z_{2}\over z_1}
-(q+1/q) 
\, + \,
{z_{m-1}\over z_{m}}
\Big)  
 \prod_{i=1}^{m-2} \prod_{j=i+2}^{m}(z_j-z_i)
\bigg)\,=\, 0.
$$
We claim that in this case we separately have
$$
Asym_m
\bigg(
 \prod_{i=1}^{m-2} \prod_{j=i+2}^{m}(z_j-z_i)
\bigg)= 0
$$

$$
a)\hskip 4pt\hskip 4pt Asym_m
\bigg(
{z_{2}\over z_1}
 \prod_{i=1}^{m-2} \prod_{j=i+2}^{m}(z_j-z_i)
\bigg)= 0
\,,\,\hskip 4pt\hskip 4pt\hskip 4pt\hskip 4pt\hskip 4pt\hskip 4pt
b)\hskip 4pt\hskip 4pt Asym_m
\bigg(
{z_{m-1}\over z_{m}} 
 \prod_{i=1}^{m-2} \prod_{j=i+2}^{m}(z_j-z_i)
\bigg)= 0.
 $$
Now \ref{eq:3.2}8 is immediate since the polynomial that
is anti-symmetrized is of degree $m-1\choose 2$.
We will complete our proof of \ref{eq:3.22} by showing 
\ref{eq:3.2}9 a). The identity in  \ref{eq:3.2}9 b) can be dealt with in an entirely analogous manner.

By collecting terms with respect to $z_1$, we have
$$
\prod_{i=1}^{m-2} \prod_{j=i+2}^m (z_j-z_i) = z_3\cdots z_m \prod_{i=2}^{m-2} \prod_{j=i+2}^m (z_j-z_i) + \sum_{k=1}^{m-2} z_1^k P_k(z_2,\dots, z_m) , 
$$
where $P_k$ is a polynomial for each $k\ge 1$. It follows that for all $k\ge 1$, $z_2/z_1 \cdot z_1^k P_k $ is a polynomial of degree ${m-1\choose 2}$ and hence they $S_m$-antisymmetrizes to $0$. For the only remaining term, we observe that
\begin{align*}
  Asym_m &\Big(z_2/z_1 \cdot z_3\cdots z_m \prod_{i=2}^{m-2} \prod_{j=i+2}^m (z_j-z_i)\Big) 
  \cr
 =&Asym_m \left({1\over z_1} z_2 z_3\cdots z_m
  Asym_{2,m} \Big(\prod_{i=2}^{m-2} \prod_{j=i+2}^m (z_j-z_i)\Big)\right)=0.
\end{align*}

Since 
$$
 Asym_{2,m} \Big(\prod_{i=2}^{m-2} \prod_{j=i+2}^m (z_j-z_i)\Big)
 \,=\, 0
$$
holds true for the same reason we have  \ref{eq:3.2}8. This completes our proof.\\

\noindent\textbf{Remark 3.1}

Convincing MAPLE to deliver zero after symmetrization in \ref{eq:3.20} is not trivial even when $m$ is as small  as $5$. We  actually succeeded in pushing the verification of \ref{eq:3.20} for all co-prime pairs $(m,n)$ with $1\le n\le m-1$
and $m\le 7$. This given, it is worth while sketching at least what we 
did for $m=5$. The cases $m=6,7$ use only more elaborate 
versions of the same ideas. 

The first step is to notice we may write
\begin{align*}
\prod_{1\le i<j\le m}\Omega[-\textstyle{z_{i}\over z_{ j}}M]
&=
\prod_{1\le i<j\le m} {(z_j-z_i)(z_j-qtz_i)\over 
(z_j-tz_i)(z_j-qz_i)}
\cr
&=
\prod_{1\le i<j\le m} {(z_j-z_i)
(z_j-qtz_i)(z_i-tz_j)(z_i-qz_j)\over 
(z_j-tz_i)(z_j-qz_i)(z_i-tz_j)(z_i-qz_j)}.
\end{align*}

Since the expression
$$
\prod_{1\le i<j\le m}(z_j-tz_i)(z_j-qz_i)(z_i-tz_j)(z_i-qz_j)
$$
is symmetric in $z_1,z_2,\dots,z_m$, it may be omitted in \ref{eq:3.20} and the
Negut identity may be also be established by proving that
$$
Asym_m
\bigg(
{\Xi_{m,n}[Z_m]
}
\Big(
\textstyle{z_{c+1}\over z_c}
-(t+q)
\, + \,
qt\textstyle{z_a\over z_{a+1}}
\Big)
\prod_{1\le i<j\le m}\hskip -.1in\hskip -.1in\,\,
(z_j-qtz_i)(z_i-tz_j)(z_i-qz_j)
\bigg)\,=\, 0.
$$
Now recall  that we may write, for ${n\over m}<1$,
$$
\Xi_{m,n}[Z_m]=
\prod_{i=2}^m
{z_i\over  z_i^{\lfloor i\textstyle{n\over m}\rfloor-\lfloor (i-1)\textstyle{n\over m}\rfloor }}\prod_{i=1}^{m-1}{1\over z_{i+1}-qtz_i}.
$$
Observing that the expression
$$
\mathcal{E}= \prod_{2\le   j\le m}\hskip -.1in\hskip -.1in\,\,
(z_1-tz_j)(z_1-qz_j)
$$
is symmetric in $z_2,z_3,\ldots z_m$, to prove \ref{eq:3.3}2 we may start by anti-symmetrizing, with respect the symmetric group 
$S_{2,m}$. The expression
$$
\mathcal{F}=\Big(
\textstyle{z_{c+1}\over z_c}
-(t+q)
\, + \,
qt\textstyle{z_a\over z_{a+1}}
\Big)
\prod_{i=2}^m
{z_i\over  z_i^{\lfloor i\textstyle{n\over m}\rfloor-\lfloor (i-1)\textstyle{n\over m}\rfloor }}
\prod_{i=1}^{m-2}
\prod_{j=i+2}^m
(z_j-qtz_i)(z_i-tz_j)(z_i-qz_j)
$$
and to prove \ref{eq:3.3}0 we now are reduced to checking  that
$$
Asym_m\bigg(\mathcal{E}\Big(Asym_{2,m} \mathcal{F}
\bigg)=0.
$$
However, to save on memory usage, we can do better than computing $Asym_{2,m} \mathcal{F}$. In fact, noticing that $\mathcal{F}$
is a Laurent polynomial we need only rewrite it in what we shall refer to as a \emph{normal form}. More precisely this amounts to
replacing  $Asym_{2,m} \mathcal{F}$ by the Laurent polynomial
$\mathcal{NFF}$ obtained by removing from $\mathcal{F}$ all the monomials with repeated exponents and then replacing each of the remaining monomials 
by the rearrangement that makes the exponents decrease, 
multiplied by the sign of the permutation that produces that 
rearrangement. Since the $S_{2,m}$ anti-symmetrization
of such a normalized monomial produces the same polynomial
yielded by the original monomial, it follows that there is no loss in replacing $Asym_{2,m}\mathcal{F}$ by $\mathcal{NFF}$ in \ref{eq:3.3}1.

It turns out that the reduction in size caused  by the  combination of these simple tricks makes MAPLE recognize that
$$
Asym_m\bigg(
\mathcal{E}  Asym_{2,m}\mathcal{NFF}
\bigg)\,=\, 0
$$
at least for $m=5$. For $m=6,7$ further partial anti-symmetrizations are necessary but the basic idea is to reduce the size as much as possible within successive anti-symmetrizations.

\section{APPENDIX}

\centerline\textbf{  The computation of the commutator $\bf D_aD_b^*-D_b^*D_a.$}

 We should mention that the identity proved here  was originally
 obtained using the Theory of Constant Terms developed in \cite{pl21}.
What we give here is a completely elementary  proof worked
out  for an audience that is unfamiliar with the above mentioned  theory.

We will adopt the following convention: for $E_i[t_1,t_2,\ldots]$ any rational functions of
the variables $t_1,t_2,\ldots$ and $P$ a symmetric polynomial, we set
$$
P^{(r_1,r_2,\ldots ,r_k)}[X]\,=\, 
P[X +\textstyle{E_1 u_1}+\textstyle{E_2 u_2}+\cdots +\textstyle{E_k u_k}]
\Big|_{u_1^{r_1}u_2^{r_2}\cdots u_k^{r_k}}.
$$
The  important property is that if
$$
Q^{(s_1)}[X]\,=\, P[X +\textstyle{E_1 u_1}]\Big|_{u_1^{s_1}}
$$
then
$$
Q^{(s_1)}[X+E_2 u_2]\Big|_{u_2^{s_2}}\,=\, P[X+E_1u_1+E_2u_2]\Big|_{u_1^{s_1}u_2^{s_2}}.
$$

For $P$ a homogeneous symmetric polynomial of degree $d$ we have
\begin{align*}
D_b^*\, P[X] 
& \,=\, \hskip 6pt  P\big[X\, - \, \, {\textstyle {\widetilde{M} \over z_2}}\, 
 \big]\, \Omega[z_2\, X\, ]\hskip 6pt \big|_{z_2^b}
\cr
& \,=\, \hskip 6pt \sum_{r_2=0}^d P^{(r_2)} 
[X]\,(\textstyle{1\over z_2})^{r_2} \sum_{u\ge 0}z_2^u \, h_u[X]\hskip 6pt \big|_{z_2^b}
\,=\, \hskip 6pt \sum_{r_2=0}^d P^{(r_2)} 
[X]\, h_{r_2+b}[X].
\end{align*}

Thus
\begin{align}
D_aD_b^*\, P[X] 
& \,=\,
\sum_{r_1,r_2=0}^d P^{(r_1,r_2)}[X](\textstyle{1\over z_1})^{r_1}
h_{r_2+b}[X+\textstyle{M\over z_1}]\Omega[-z_1\, X\, ]\hskip 6pt \Big|_{z_1^a}
\cr
& \,=\,
\sum_{r_1,r_2=0}^d P^{(r_1,r_2)}[X](\textstyle{1\over z_1})^{r_1}
\sum_{s=0}^{r_2+b}h_{r_2+b-s}[X](\textstyle{1\over z_1})^s h_s[M]
\sum_{u\ge 0}z_1^uh_u[-X]\hskip 6pt \Big|_{z_1^a}
\cr
& \,=\,
\sum_{r_1,r_2=0}^d P^{(r_1,r_2)}[X] 
\sum_{s=0}^{r_2+b}h_{r_2+b-s}[X]  h_s[M]
 h_{r_1+s+a}[-X].
\label{eq:A.1}
\end{align}

Making the summation parameter change $u=r_1+s+a$ gives $s=u-r_1-a$
and the range
$$
 r_1 +a\le u\le  r_1+r_2+a+b
$$
so \ref{eq:A.1} becomes
$$
D_aD_b^*\, P[X]\,=\, \sum_{r_1,r_2=0}^d P^{(r_1,r_2)}[X] 
\sum_{u=r_1 +a}^{r_1+r_2+a+b}h_{r_1+r_2+a+b -u }[X] h_{u}[-X] 
h_{u-r_1-a}[M].
$$
Similarly
\begin{align*}
D_a \, P[X] 
& \,=\, \hskip 6pt  P\big[X+ {\textstyle { M \over z_1}}\, 
 \big]\, \Omega[-z_1 X  ]\hskip 6pt \Big|_{z_1^a}
\cr
& \,=\, \hskip 6pt \sum_{r_1=0}^d P^{(r_1)} 
[X]\,(\textstyle{1\over z_1})^{r_1} \sum_{u\ge 0}z_1^u \, h_u[-X]\hskip 6pt \big|_{z_1^a}
\,=\, \hskip 6pt \sum_{r_1=0}^d P^{(r_1)} 
[X]\, h_{r_1+a}[-X].
\end{align*}

Thus
\begin{align}
D_b^*D_a\, P[X] 
& \,=\,
\sum_{r_1,r_2=0}^d P^{(r_1,r_2)}[X](\textstyle{1\over z_2})^{r_2}
h_{r_1+a}[-X+\textstyle{\widetilde{M}\over z_2}]\Omega[z_2\, X\, ]\hskip 6pt \Big|_{z_2^b}
\cr
& \,=\,
\sum_{r_1,r_2=0}^d P^{(r_1,r_2)}[X](\textstyle{1\over z_2})^{r_2}
\sum_{s=0}^{r_1+a}h_{r_1+a-s}[-X](\textstyle{1\over z_2})^s h_s[\widetilde{M}]
\sum_{u\ge 0}z_2^uh_u[X]\hskip 6pt \Big|_{z_2^b}
\cr
& \,=\,
\sum_{r_1,r_2=0}^d P^{(r_1,r_2)}[X] 
\sum_{s=0}^{r_1+a}h_{r_1+a-s}[-X]  h_s[\widetilde{M}]
 h_{r_2+s+b}[X].
\label{eq:A.2}
\end{align}

Making the summation parameter change $u=r_1+a-s$ gives $s= r_1+a-u$
and the range
$$
0  \le u \le  r_1+a 
$$
ao \ref{eq:A.2} becomes
$$
D_b^*D_a\, P[X]\,=\,
\sum_{r_1,r_2=0}^d P^{(r_1,r_2)}[X] 
\sum_{u=0}^{r_1+a}h_{u}[-X] 
 h_{r_1+r_2+ a+b -u }[X] h_{r_1+a-u  }[\widetilde{M}].
$$
Now recall that 
$$
h_m[M]= M\textstyle{1-t^mq^m\over 1-  tq}
\hskip 4pt \hskip 4pt\hskip 4pt\hskip 4pt\hskip 4pt \hbox{and} \hskip 4pt \hskip 4pt\hskip 4pt\hskip 4pt\hskip 4pt
h_m[\widetilde{M}]= M\textstyle{1-t^mq^m\over 1-  tq}\textstyle{1\over t^mq^m}
= -M\textstyle{1-1/t^mq^m\over 1-  tq}.
$$
We thus get
$$
D_aD_b^*\, P[X]\,=\, \textstyle{M\over 1-tq}\sum_{r_1,r_2=0}^d P^{(r_1,r_2)}[X] 
\sum_{u=r_1 +a}^{r_1+r_2+a+b}h_{r_1+r_2+a+b -u }[X] h_{u}[-X] 
(1-(tq)^{{u-r_1-a}})
$$
and
$$
D_b^*D_a\, P[X]\,=\,
-\textstyle{M\over 1-tq}\sum_{r_1,r_2=0}^d P^{(r_1,r_2)}[X] 
\sum_{u=0}^{r_1+a}h_{u}[-X] 
 h_{r_1+r_2+ a+b -u }[X](1-(tq)^{{u-r_1-a}}).
$$
Hence
\begin{align*}
(D_aD_b^*-D_b^*D_a)\, P[X]
&\,=\, \textstyle{M\over 1-tq} \displaystyle \sum_{r_1,r_2=0}^d P^{(r_1,r_2)}[X] 
\hskip -.15in \sum_{u=0}^{r_1+r_2+a+b} \hskip -.15in
h_{r_1+r_2+a+b -u }[X] h_{u}[-X] 
(1-(tq)^{{u-r_1-a}})
\cr
&\,=\, \textstyle{M\over 1-tq} \displaystyle\sum_{r_1,r_2=0}^d P^{(r_1,r_2)}[X] 
 h_{r_1+r_2+a+b }[X-X] 
\cr
&\hskip .25in
\, - \, \textstyle{M\over 1-tq} \displaystyle \sum_{r_1,r_2=0}^d P^{(r_1,r_2)}[X] 
(tq)^{{ -r_1-a}}
\hskip -.15in \sum_{u=0}^{r_1+r_2+a+b} \hskip -.15in
h_{r_1+r_2+a+b -u }[X] h_{u}[-tqX]  
\cr
&\,=\, \textstyle{M\over 1-tq}\sum_{r_1,r_2=0}^d P^{(r_1,r_2)}[X] 
 h_{r_1+r_2+a+b }[X-X] 
\cr
&\hskip .25in
\, - \, \textstyle{M\over 1-tq} \displaystyle \sum_{r_1,r_2=0}^d P^{(r_1,r_2)}[X] 
(tq)^{{ -r_1-a}} h_{r_1+r_2+a+b   }[X(1-tq)].
\end{align*}

But now note that we may write
\begin{align*}
\sum_{r_1,r_2=0}^d  P^{(r_1,r_2)}[X] &
(tq)^{{ -r_1-a}} h_{r_1+r_2+a+b   }[X(1-tq)] \cr
&=
\textstyle{1\over (tq)^a } \sum_{r_1,r_2=0}^d  P[X+M u_1-\widetilde{M} u_2]\Big|_{u_1^{r_1} u_2^{r_2}}
(\textstyle{1\over tqz})^{r_1} (\textstyle{1\over z})^{r_2}  \Omega[zX(1-tq)]\Big|_{z^{a+b}}
\cr
&=
\textstyle{1\over (tq)^a }      P[X+\textstyle{M\over qtz}  -\textstyle{\widetilde{M}\over z } ] \hskip 4pt
   \Omega[zX(1-tq)]\Big|_{z^{a+b}}
\cr
&=
\textstyle{1\over (tq)^a }      P[X] \hskip 4pt
   \Omega[zX(1-tq)]\Big|_{z^{a+b}}.
\end{align*}

This proves
\begin{equation}
(D_aD_b^*-D_b^*D_a)\, P[X]\,=\, 
\textstyle{M\over 1-tq}\begin{cases}
{-1\over (qt)^a}h_{a+b}\big[X(1-qt)\big]
P[X] & \hbox{if }a+b>0,
\cr
\big(1-{1\over (qt)^a}\big) P[X]
& \hbox{if } a+b=0,
\cr
\sum_{r_1+r_2=-(a+b)}P^{(r_1,r_2)}[X]
& \hbox{if } a+b<0.
\end{cases}
\label{eq:A.3}
\end{equation}

\nocite{*}
\bibliographystyle{plain} 
\bibliography{SymCompos2014}

\end{document}